\documentclass{article}%
\usepackage{amsmath}
\usepackage{amsfonts}
\usepackage{amssymb}
\usepackage{graphicx}
\usepackage{tikz,pgfplots}
\usepackage{epstopdf}
\usepackage{color}%
\providecommand{\U}[1]{\protect\rule{.1in}{.1in}}
\providecommand{\U}[1]{\protect\rule{.1in}{.1in}}
\newtheorem{theorem}{Theorem}

\newtheorem{corollary}[theorem]{Corollary}
\newtheorem{definition}[theorem]{Definition}

\newtheorem{lemma}[theorem]{Lemma}

\newtheorem{remark}[theorem]{Remark}

\newenvironment{proof}[1][Proof]{\noindent\textbf{#1.} }{\ \rule{0.5em}{0.5em}}
\begin{document}

\title{A Fully Discrete Galerkin Method for Abel-type Integral Equations}
\author{Urs V\"{o}geli\thanks{Institut f\"{u}r Mathematik, Universit\"{a}t Z\"{u}rich,
Winterthurerstrasse 190, CH-8057 Z\"{u}rich, Switzerland, e-mail:
\texttt{voegeli.urs@gmail.com}}
\and Khadijeh Nedaiasl\thanks{Institute for Advanced Studies in Basic Sciences,
Zanjan, Iran, e-mail: \texttt{nedaiasl@iasbs.ac.ir}}
\and Stefan A. Sauter\thanks{Institut f\"{u}r Mathematik, Universit\"{a}t
Z\"{u}rich, Winterthurerstrasse 190, CH-8057 Z\"{u}rich, Switzerland, e-mail:
\texttt{stas@math.uzh.ch}}}
\maketitle

\begin{abstract}
In this paper, we present a Galerkin method for Abel-type integral equation
with a general class of kernel. Stability and quasi-optimal convergence
estimates are derived in fractional-order Sobolev norms. The fully-discrete
Galerkin method is defined by employing simple tensor-Gauss quadrature. We
develop a corresponding perturbation analysis which allows to keep the number
of quadrature points small. Numerical experiments have been performed which
illustrate the sharpness of the theoretical estimates and the sensitivity of
the solution with respect to some parameters in the equation.

\end{abstract}

\vspace{1em} \noindent\textbf{Keywords:} Abel's integral equation, Galerkin
method, tensor-Gauss quadrature

\vspace{1em} \noindent\textbf{Mathematics Subject Classification (2000):}
45E10, 65R20, 65D32

\section{Introduction}

A variety of practical physical models, e.g., in thermal tomography,
spectroscopy, astrophysics, can be modelled by Abel-type integral equations
provided the problem enjoys symmetries which allow to reduce the equation to a
one-dimensional equation (cf. \cite{gorenflo91}, \cite{TTot}). In this paper,
we present a fully discrete Galerkin method for the numerical solution of
Abel-type integral equations. The theory and discretization of Abel's integral
equation have been investigated by many authors and a variety of methods are
proposed which include, e.g., product integration method \cite{cameron84} and
approximation by implicit interpolation \cite{brunner04alt}. In
\cite{bieniasz08} the solution is approximated by means of an adaptive Huber
method which is a kind of product integration method. A Nystr\"{o}m-type
method which is based on the trapezoidal rule is analyzed for the Abel
integral equation in \cite{eggermont1981new}. Recently, the composite
trapezoidal method is applied for weakly singular Volterra integral equations
of the first kind \cite{plato12}. Furthermore, the piecewise polynomial
discontinuous Galerkin approximation of a first-kind Volterra integral
equation of convolution kernel type for a smooth kernel is studied in
\cite{brunner09}. The stability and robustness of the collocation method for
Abel integral equation have been discussed in \cite{eggermont84}. Less
research has been devoted to Galerkin discretizations and to the development
of a stability and convergence theory in energy spaces which are
fractional-order Sobolev spaces. {In our paper, we will propose a variational
formulation of \textit{generalized} Abel's integral equation and prove
continuity and coercivity in the natural energy spaces which are, for these
applications, fractional order Sobolev spaces. This allows to employ the
classical Lax-Milgram theory to derive well-posedness for the continuous
problem as well as for its Galerkin discretization. Then, C\'{e}a's lemmas
implies quasi-optimality of the proposed method.} For general Abel-type kernel
functions, the arising double integration for computing the coefficients in
the Galerkin system matrix, in general, cannot be evaluated analytically and
numerical quadrature has to be employed. In our paper, we present a fully
discrete Galerkin method via numerical quadrature and present a convergence
analysis in energy norms for the exact Galerkin method as well as for the
fully discrete version by employing a perturbation theory, {which allows for
numerical quadrature and a quadrature error analysis for the regular and
nearly-singular integrals. The errors arising from the approximation of the
\textit{singular} integrals can be treated within the same perturbation theory
while the derivation of these local quadrature errors require certain
regularization techniques (cf. }\cite[Chap. 5]{sauter11}). We omit the details
here in order not to overload this paper.

This perturbation theory will also allow to estimate the accuracy of
\textit{fast methods}, e.g., multipole methods for the sparse approximation of
the non-local integral operator (see, e.g., \cite{sauter11}), although this is
also beyond of the scope of this paper.

The paper is structured as follows. In Section \ref{Sec2}, we introduce the
class of Abel-type integral equations which we will investigate and present
its variational formulation in appropriate fractional-order Sobolev spaces. We
will show well-posedness of the problem by proving continuity and ellipticity
of the arising sesquilinear form and employ the Lax-Milgram lemma. We employ a
lemma from \cite{eggermont88} for the ellipticity of the original Abel
integral equation and then derive the result for generalized types of integral
kernels by a perturbation argument.

In Section \ref{Sec3}, we present the Galerkin discretization and introduce
the one-dimensional polynomial finite element spaces as test and trial spaces.
Section \ref{Sec4} is devoted to the stability and convergence analysis of the
Galerkin method. Then, we introduce a simple numerical quadrature method based
on tensor-Gauss formulae and develop a perturbation analysis which allows to
keep the number of quadrature points small.

In Section \ref{SecNumexp}, we report on numerical experiments which
illustrate that the theoretically derived convergence estimates are sharp and
also the estimates for the number of quadrature nodes is close to optimal.
Finally, we vary the parameter $\alpha\in\left(  0,1\right)  $ in Abel's
integral equation which determines the strength of singularity of the kernel
function to study numerical its effect on the solution.

In the concluding Section \ref{SecConcl} we briefly summarize our results.

\section{Abel-type Integral Equations \label{Sec2}}

\subsection{Setting}

In this paper, we consider Abel-type integral equations of the following
abstract form: Let $\Omega=\left(  0,1\right)  $. For a given function
$g:\Omega\rightarrow\mathbb{R}$, parameter $\alpha\in\left(  0,1\right)  $,
and kernel function $K:\Omega\times\Omega\rightarrow\mathbb{R}$ we are seeking
the solution $f$ of the equation%
\begin{equation}
A_{K,\alpha}(f)(x):=\frac{1}{\Gamma(\alpha)}\int_{0}^{x}{(x-y)^{\alpha
-1}K(x,y)f(y)dy}=g(x),\quad\text{for all }x\in\Omega, \label{inteq}%
\end{equation}
where $A_{K,\alpha}$ is the Abel-type integral operator and $\Gamma(\alpha)$
denotes Euler's Gamma function. To formulate the appropriate function spaces
so that (\ref{inteq}) is well posed we first recall the definition of
fractional-order Sobolev spaces. As usual the standard Lebesgue spaces are
denoted by $L^{p}\left(  \Omega\right)  $ and their norm by $\left\Vert
\cdot\right\Vert _{L^{p}\left(  \Omega\right)  }$. For $p=2$, the scalar
product is denoted by $\left(  u,v\right)  =\int_{\Omega}u\overline{v}$ and
the norm by $\left\Vert \cdot\right\Vert =\left(  \cdot,\cdot\right)  ^{1/2}$.

For $\ell\in\mathbb{R}$, let $\left\lfloor \ell\right\rfloor $ denote the
largest integer for which $\left\lfloor \ell\right\rfloor \leqslant\ell$ and
define $\lambda\in\left[  0,1\right[  $ by $\ell=\left\lfloor \ell
\right\rfloor +\lambda$. For $\ell\in\mathbb{R}_{>0}\backslash\mathbb{N}$, we
introduce the scalar product%
\begin{align}
\left(  \varphi,\psi\right)  _{H^{\ell}(\Omega)} &  :=\sum_{\alpha
\leqslant\left\lfloor \ell\right\rfloor }{\left(  D^{\alpha}\varphi,D^{\alpha
}\psi\right)  }\label{lscalar}\\
&  +\sum_{\alpha=\left\lfloor \ell\right\rfloor }{\int_{\Omega\times\Omega
}{\frac{\left(  D^{\alpha}\varphi(x)-D^{\alpha}\varphi(y)\right)
\overline{\left(  D^{\alpha}\psi(x)-D^{\alpha}\psi(y)\right)  }}%
{|x-y|^{1+2\lambda}}dxdy}}\nonumber
\end{align}
and the norm $\Vert\varphi\Vert_{H^{\ell}(\Omega)}:=\left(  \varphi
,\varphi\right)  _{H^{\ell}(\Omega)}^{1/2}$. For $\ell\in\mathbb{N}$, the
second term in (\ref{lscalar}) is skipped. Then the Sobolev space $H^{\ell
}\left(  \Omega\right)  $ is given by
\[
H^{\ell}\left(  \Omega\right)  :=\left\{  u\in L^{2}(\Omega)\mid\forall0\leq
m\leq\left\lfloor \ell\right\rfloor \quad u^{\left(  m\right)  }\in
L^{2}\left(  \Omega\right)  \text{\quad and\quad}\Vert u\Vert_{H^{\ell}%
(\Omega)}<\infty\right\}  .
\]
For our application the differentiation order $\ell$ will satisfy $\ell<1/2$.
In this case, the dual space of $H^{\ell}(\Omega)$ is denoted by $H^{-\ell
}(\Omega)$ and is equipped with the {norm}
\begin{equation}
\Vert u\Vert_{H^{-\ell}(\Omega)}:=\sup\limits_{v\in H^{\ell}\left(
\Omega\right)  \backslash\left\{  0\right\}  }\frac{(u,v)}{\Vert
v\Vert_{H^{\ell}(\Omega)}},
\end{equation}
where $(\cdot,\cdot)$ denotes the continuous extension of the $L^{2}$ scalar
product to the anti-dual pairing $\langle\cdot,\overline{\cdot}\rangle$ in
$H^{-\ell}(\Omega)\times H^{\ell}(\Omega)$.

\subsection{Variational Formulation of Abel's Integral Equation}

{To state the variational form of the Abel's integral equation, we consider
(\ref{inteq}) as an operator equation in $H^{\alpha/2}\left(  \Omega\right)  $
which will be justified by Theorems \ref{2.3} and \ref{main1}: } For given
$g\in H^{\alpha/2}\left(  \Omega\right)  $, we are seeking $f\in H^{-\alpha
/2}\left(  \Omega\right)  $ such that%
\begin{equation}
A_{K,\alpha}f=g. \label{gal1}%
\end{equation}
We multiply (\ref{gal1}) by test functions $\varphi\in H^{-\alpha/2}\left(
\Omega\right)  $ and integrate over $\Omega$ to get the variational problem:

Find $f\in H^{-\alpha/2}(\Omega)$ such that%
\begin{equation}
a_{K,\alpha}\left(  f,\varphi\right)  :=\left(  A_{K,\alpha}f,\varphi\right)
=(g,\varphi)=:G\left(  \varphi\right)  ,\quad\forall\varphi\in H^{-\alpha
/2}(\Omega). \label{var1}%
\end{equation}

\subsection{Well-posedness of the Variational Problem}

In this section we prove the well-posedness of (\ref{var1}) under certain
assumptions on the kernel $K$: In Section \ref{SubSecEll1}, we recall the
ellipticity of the sesquilinear form $a_{K,\alpha}$ for $K=1$ and extend this
result to more general kernels in Section \ref{EllContgenK}, i.e., there
exists a constant $\gamma>0$ such that%
\begin{equation}
\operatorname{Re}a_{K,\alpha}(u,u)\geqslant\gamma\Vert u\Vert_{H^{-\alpha
/2}\left(  \Omega\right)  }^{2}\quad\forall u\in H^{-\alpha/2}\left(
\Omega\right)  . \label{helliptic}%
\end{equation}
In Section \ref{SecConti} we prove the continuity of $a_{K,\alpha}$ for $K=1$
based on results in \cite{gorenflo99} and extend this results for a larger
class of kernels in Section \ref{EllContgenK}. Then, well-posedness of
(\ref{var1}) follows from the Lax-Milgram lemma.


\subsubsection{Ellipticity \label{SubSecEll1}}

We first recall the $H^{-\alpha/2}$-ellipticity for the case $K=1$ and
generalize this result to a much larger class of kernel functions $K\neq1$ by
a perturbation argument.

\begin{theorem}
We consider $A_{K,\alpha}$ in (\ref{inteq}) for $K=1$ and $\alpha\in\left(
0,1\right)  $. The operator $A_{1,\alpha}$ is $H^{-\alpha/2}(\Omega
)$-elliptic:
\begin{equation}
\operatorname{Re}\left(  A_{1,\alpha}f,f\right)  \geqslant\gamma\Vert
f\Vert_{H^{-\alpha/2}\left(  \Omega\right)  }^{2}\quad\forall f\in
H^{-\alpha/2}\left(  \Omega\right)  \quad\text{with }\gamma:=\cos\left(
\frac{\pi\alpha}{2}\right)  . \label{ellineq}%
\end{equation}
\label{coercivity}
\end{theorem}

\proof
In \cite[Lem. 2.3]{eggermont88}, it is shown that the inequality
(\ref{ellineq}) holds for $f\in L^{2}(\Omega)$. Since $L^{2}(\Omega)$ is dense
in $H^{-\alpha/2}(\Omega)$ (cf. \cite[Prop. 2.4.2]{sauter11}) this also holds
for all $f\in H^{-\alpha/2}(\Omega)$.\endproof

\subsubsection{Continuity \label{SecConti}}

In this section we prove the continuity of $a_{1,\alpha}$ as in (\ref{var1}).
For the continuity we will follow the proof in \cite{gorenflo99} and start
with some preliminaries. Let
\begin{equation}
\mu_{n}=\left(  n-\frac{1}{2}\right)  \pi,\quad\phi_{n}(t)=\sqrt{2}\cos\mu
_{n}t,\text{ for }0\leqslant t\leqslant1,n\in\mathbb{N}_{\geqslant1}.
\end{equation}

\begin{remark}
Note that $\left(  \mu_{n}^{2},\phi_{n}\right)  _{n\geq1}$ is the set of all
eigenpairs of the boundary-value problem,
\begin{equation}%
\begin{split}
(1)\qquad &  (D^{2}u)(t)=-\lambda u(t),\quad0<t<1,\\
(2)\qquad &  \frac{du}{dt}(0)=u(1)=0.
\end{split}
\label{eq:19}%
\end{equation}
It is known that $(\phi_{n})_{n\in\mathbb{N}_{\geqslant1}}$ is complete in
$L^{2}(\Omega)$ and the eigenfunctions $\left(  \phi_{n}\right)  _{n\geq1}$
form an orthonormal basis in $L^{2}\left(  \Omega\right)  $ (cf. \cite[Eq.
(2.1)]{gorenflo99}).
\end{remark}

For $\beta\in\mathbb{R}$, we may introduce as in \cite[Eq. (2.2)]{gorenflo99}
a Hilbert scale $X_{\beta}$ by employing this basis. In $\operatorname*{span}%
\left\{  \phi_{n}:n\geq1\right\}  $ we define the scalar product and norm by
\begin{equation}
(u,v)_{X_{\beta}}:=\sum_{n=1}^{\infty}{\mu_{n}^{2\beta}(u,\phi_{n}%
)\overline{(v,\phi_{n})}},\quad\Vert u\Vert_{X_{\beta}}:=(u,u)_{X_{\beta}%
}^{1/2},\quad\forall u,v\in\operatorname*{span}\left\{  \phi_{n}%
:n\geq1\right\}  . \label{scalenorm}%
\end{equation}

We introduce an operator $S:L^{2}(\Omega)\rightarrow L^{2}(\Omega)$ by
\begin{equation}
(Su)(t)=\sum_{n=1}^{\infty}{\mu_{n}^{-1}(u,\phi_{n})\phi_{n}(t)},\quad u\in
L^{2}(\Omega).
\end{equation}
For the proofs of the following two lemmata, we refer to \cite[Lem. 7,
8]{gorenflo99}.

\begin{lemma}
The fractional power $S^{\alpha}$ of $S$ is given by
\begin{equation}
S^{\alpha}u= \sum_{n=1}^{\infty}{\mu_{n}^{-\alpha}(u,\phi_{n})\phi_{n}},\quad
u\in L^{2}(\Omega), \,\alpha>0.
\end{equation}
\label{salpha}
\end{lemma}

The next lemma gives further insights on the space $X_{\beta}$.

\begin{lemma}
\hfill

\begin{itemize}
\item $X_{\beta}=H^{\beta}(\Omega)$, $0\leqslant\beta<\frac{1}{2}$,

\item $X_{1/2}=\{u\in H^{1/2}(\Omega):\int_{0}^{1}{(1-t)^{-1}|u(t)|^{2}%
dt}<\infty\}$,

\item $X_{\beta}=\{u\in H^{\beta}(\Omega): u(1)=0\}$, $\frac{1}{2}%
<\beta\leqslant1.$
\end{itemize}

Moreover there exists a constant $C=C(\beta)>0$ such that
\begin{equation}
C^{-1}\Vert u\Vert_{X_{\beta}}\leqslant\Vert u\Vert_{H^{\beta}\left(
\Omega\right)  }\leqslant C\Vert u\Vert_{X_{\beta}},\quad u\in X_{\beta},
\end{equation}
if $\beta\in\lbrack0,1]$ and $\beta\neq1/2$. \label{XX}
\end{lemma}

This lemma is used in the proof of the continuity for $A_{K,\alpha}$.

\begin{theorem}
\label{2.3} Let the kernel function $K=K(x,y)$ satisfy the following conditions:

\begin{itemize}
\item $K$ is continuous on $D=\{(x,y)\in\mathbb{R}^{2}:0\leqslant y\leqslant
x\leqslant1\}$,

\item $K(x,x)=1$ for $0\leqslant x\leqslant1$,

\item there exists a decreasing function $k\in L^{2}(\Omega)$ such that
$\left\vert \frac{\partial K}{\partial{y}}(x,y)\right\vert \leqslant k(y)$ for
$0<y\leqslant x\leqslant1$. \label{kcond}
\end{itemize}

Abel's integral operator $A_{K,\alpha}$ is a continuous mapping
from$\ H^{-\alpha/2}(\Omega)$ into $H^{\alpha/2}(\Omega)$. \label{continuity}
\end{theorem}

%

\proof
We show that the operator from $H^{-\alpha/2}\left(  \Omega\right)  $ to
$H^{\alpha/2}\left(  \Omega\right)  $ is bounded so we can conclude its
continuity. From \cite[Thm. 1]{gorenflo99}, we know that there exists a
constant $C=C(\alpha)$ for $\alpha\in(0,1)$ such that
\begin{equation}
C^{-1}\Vert u\Vert_{X_{-\alpha}}\leqslant\Vert A_{K,\alpha}u\Vert\leqslant
C\Vert u\Vert_{X_{-\alpha}},\quad u\in L^{2}(\Omega). \label{aju}%
\end{equation}
Next, we investigate the self-adjoint compact operator $(A_{K,\alpha}^{\ast
}A_{K,\alpha})^{1/2}:\newline L^{2}(\Omega)\rightarrow L^{2}(\Omega)$ with its
eigenvalues $(s_{n}(A_{K,\alpha}))_{n\in\mathbb{N}_{\geqslant1}}$. From
\cite[Thm. 2]{gorenflo99}, we conclude that there exists a constant
$C=C(\alpha)$ such that%
\begin{equation}
C^{-1}n^{-\alpha}\leqslant s_{n}(A_{K,\alpha})\leqslant Cn^{-\alpha}.
\label{sing}%
\end{equation}
From (\ref{scalenorm}) we obtain ${{(\phi_{n},\phi_{m})_{X_{\alpha/2}}=\delta
}}_{n,m}\mu_{n}^{2}={\delta}_{n,m}\left\Vert \phi_{n}\right\Vert
{{_{X_{\alpha/2}}^{2}}}$. Hence%
\begin{align*}
\Vert A_{K,\alpha}u\Vert_{X_{\alpha/2}}^{2}  &  =(A_{K,\alpha}u,A_{K,\alpha
}u)_{X_{\alpha/2}}=\sum_{n=1}^{\infty}{\sum_{m=1}^{\infty}{u_{n}%
\overline{u_{m}}(A_{K,\alpha}\phi_{n},A_{K,\alpha}\phi_{m})_{X_{\alpha/2}}}}\\
&  =\sum_{n=1}^{\infty}{\sum_{m=1}^{\infty}{u_{n}\overline{u_{m}}%
s_{n}(A_{K,\alpha})s_{m}(A_{K,\alpha})(\phi_{n},\phi_{m})_{X_{\alpha/2}}}}\\
&  =\sum_{n=1}^{\infty}{|u_{n}|^{2}s_{n}(A_{K,\alpha})^{2}\Vert\phi_{n}%
\Vert_{X_{\alpha/2}}^{2}}.
\end{align*}
From (\ref{scalenorm}) we deduce $\left\Vert \phi_{n}\right\Vert
_{X_{\alpha/2}}^{2}=\mu_{n}^{\alpha}$ so that%
\[
\Vert A_{K,\alpha}u\Vert_{X_{\alpha/2}}^{2}=\sum_{n=1}^{\infty}{|u_{n}%
|^{2}s_{n}(A_{K,\alpha})^{2}\mu_{n}^{\alpha}}.
\]
We combine this with the estimate of $s_{n}$ as in (\ref{sing}) and the
definition of $\mu_{n}$ to obtain%
\begin{align*}
\Vert A_{K,\alpha}u\Vert_{X_{\alpha/2}}^{2}  &  \leqslant C\sum_{n=1}^{\infty
}{|u_{n}|^{2}n^{-2\alpha}\mu_{n}^{\alpha}}=C\sum_{n=1}^{\infty}{|u_{n}%
|^{2}n^{-2\alpha}(\pi(n-1/2))^{\alpha}}\\
&  \leqslant C\pi^{2\alpha}\sum_{n=1}^{\infty}{|u_{n}|^{2}(\pi
(n-1/2))^{-\alpha}}=C\pi^{2\alpha}\sum_{n=1}^{\infty}{|u_{n}|^{2}\mu
_{n}^{-\alpha}}.
\end{align*}
The last ingredient for the proof is%
\begin{align*}
\Vert u\Vert_{X_{-\alpha/2}}^{2}  &  =\Vert S^{\alpha/2}u\Vert^{2}=\Vert
\sum_{n=1}^{\infty}{\mu_{n}^{-\alpha/2}(\sum_{m=1}^{\infty}{u_{m}\phi_{m}%
},\phi_{n})\phi_{n}}\Vert^{2}\\
&  =\Vert\sum_{n=1}^{\infty}{\mu_{n}^{-\alpha/2}u_{n}\phi_{n}}\Vert
^{2}=\left(  \sum_{n=1}^{\infty}{\mu_{n}^{-\alpha/2}u_{n}\phi_{n}},\sum
_{m=1}^{\infty}{\mu_{m}^{-\alpha/2}u_{m}\phi_{m}}\right) \\
&  =\sum_{n=1}^{\infty}{\mu_{n}^{-\alpha}|u_{n}|^{2}(\phi_{n},\phi_{n})}%
=\sum_{n=1}^{\infty}{\mu_{n}^{-\alpha}|u_{n}|^{2}}.
\end{align*}
The combination of these relations leads to%
\[
\Vert A_{K,\alpha}u\Vert_{X_{\alpha/2}}^{2}\leqslant C\pi^{2\alpha}\Vert
u\Vert_{X_{-\alpha/2}}^{2}.
\]
From \cite[Chap. 3]{gorenflo99} it follows%
\[
H_{0}^{\beta}(\Omega)=H^{\beta}(\Omega)\text{ and }X_{-\beta}=H^{-\beta
}(\Omega)\qquad\forall\,0\leqslant\beta<\frac{1}{2}.
\]
Since $\alpha/2\in(0,\frac{1}{2})$ we conclude from Lemma \ref{XX} that%
\begin{equation}
\Vert A_{K,\alpha}u\Vert_{H^{\alpha/2}\left(  \Omega\right)  }\leqslant
C_{c}\Vert u\Vert_{H^{-\alpha/2}\left(  \Omega\right)  }\quad\forall u\in
H^{-\alpha/2}\left(  \Omega\right)  \label{continuityestimate}%
\end{equation}
holds.%
\endproof

\subsection{Ellipticity and Continuity for $K\neq1$ \label{EllContgenK}}

The following lemma is needed for the proof of the main result of this section
and we refer for a proof to \cite{brezis}. For $0\leq s\leq\infty$ and $1\leq
p\leq\infty$, let $W^{s,p}\left(  \Omega\right)  $ denote the usual Sobolev
space as defined, e.g., in \cite{Adams}, equipped with the norm $\left\Vert
\cdot\right\Vert _{W^{s,p}\left(  \Omega\right)  }$.

\begin{lemma}
\label{brezis}Let $1<p<\infty$, $0<s<\infty$, $0<q<\infty$, $0<\theta<1$,
$1<t<\infty$ be such that
\begin{equation}
\dfrac{1}{q}+\dfrac{\theta}{t}=\dfrac{1}{p}. \label{cond}%
\end{equation}
For $f\in W^{s,t}(\Omega)\cap L^{\infty}(\Omega)$, $g\in W^{\theta s,p}%
(\Omega)\cap L^{r}(\Omega)$, we have $fg\in W^{\theta s,p}(\Omega)$ and
\begin{equation}
\rVert fg\rVert_{W^{\theta s,p}}\leq C\left(  \rVert f\rVert_{L^{\infty}%
}\rVert g\rVert_{W^{\theta s,p}}+\rVert g\rVert_{L^{q}}\rVert f\rVert
_{W^{s,t}}^{\theta}\rVert f\rVert_{L^{\infty}}^{1-\theta}\right)  .
\label{ineq1}%
\end{equation}

\end{lemma}

For $s\in\left[  0,1\right]  $, we introduce the intervals
\[
I_{s}:=\left\{
\begin{array}
[c]{ll}%
\left[  1,\dfrac{2}{1-2s}\right]  & 0\leq s<1/2,\\
\left[  1,\infty\right[  & s=1/2,\\
\left[  1,\infty\right]  & 1/2<s\leq1
\end{array}
\right.
\]
which are relevant for the embedding properties of Sobolev spaces for the
interval $\Omega=\left(  0,1\right)  $:%
\begin{equation}
H^{s}\left(  \Omega\right)  \hookrightarrow L^{q}(\Omega)\qquad\forall
s\in\left[  0,1\right]  \quad\forall q\in I_{s}. \label{Embedding}%
\end{equation}

We apply Lemma \ref{brezis} to prove norm equivalences for products of
functions in Sobolev spaces. For this we first define a class of multipliers.

\begin{definition}
Let $s\in\left[  -1,1\right]  $ and $\Omega=\left(  0,1\right)  $. A function
$\kappa:\Omega\rightarrow\mathbb{R}$ belongs to the \emph{multiplier class}
$\mathcal{M}\left(  s\right)  $ if $\kappa\in L^{\infty}(\Omega)$ and
$\underset{x\in\Omega}{\operatorname*{ess}\inf}\kappa(x)=:\kappa_{\min}>0$ and%
\[
\kappa\text{ and }\frac{1}{\kappa}\in\left\{
\begin{array}
[c]{ll}%
H^{1}\left(  \Omega\right)  & \left\vert s\right\vert =1,\\
W^{\frac{s}{1-\varepsilon},2-\frac{2\varepsilon}{2\varepsilon-1}}\left(
\Omega\right)  & 0<\left\vert s\right\vert <1\quad\text{for some }%
\varepsilon\in\left\{
\begin{array}
[c]{ll}%
\left(  0,1/2\right)  & \text{if }1/2\leq\left\vert s\right\vert <1,\\
\left[  \frac{1}{2}-\left\vert s\right\vert ,1/2\right)  & \text{if
}0<\left\vert s\right\vert <1/2.
\end{array}
\right.
\end{array}
\right.
\]

\end{definition}

\begin{lemma}
\label{spos} Let $s\in\left[  -1,1\right]  $ and assume that $\kappa
\in\mathcal{M}\left(  s\right)  $. Then there exist constants $c_{\kappa,s}$
and $C_{\kappa,s}$ such that
\begin{equation}
c_{\kappa,s}\rVert f\rVert_{{H}^{s}\left(  \Omega\right)  }\leq\rVert\kappa
f\rVert_{H^{s}\left(  \Omega\right)  }\leq C_{\kappa,s}\rVert f\rVert
_{H^{s}\left(  \Omega\right)  },\quad\forall f\in H^{s}(\Omega). \label{29}%
\end{equation}

\end{lemma}%

\proof
For $s=0$, the inequality (\ref{29}) is easily obtained by H\"{o}lder's
inequality%
\[
\kappa_{\min}\left\Vert f\right\Vert \leq\left\Vert \kappa f\right\Vert
\leq\left\Vert \kappa\right\Vert _{L^{\infty}\left(  \Omega\right)
}\left\Vert f\right\Vert .
\]
For $s=1$ we employ Leibniz' rule%
\[
\left\Vert \kappa f\right\Vert _{H^{1}\left(  \Omega\right)  }\leq\left\Vert
\kappa^{\prime}f\right\Vert +\left\Vert \kappa f^{\prime}\right\Vert
\leq\left\Vert \kappa^{\prime}f\right\Vert +\left\Vert \kappa\right\Vert
_{L^{\infty}\left(  \Omega\right)  }\left\Vert f^{\prime}\right\Vert .
\]
A H\"{o}lder's inequality for $q^{-1}+\left(  q^{\prime}\right)  ^{-1}=1$
leads to%
\[
\left\Vert \kappa^{\prime}f\right\Vert \leq\left\Vert \kappa^{\prime
}\right\Vert _{L^{2q^{\prime}}\left(  \Omega\right)  }\left\Vert f\right\Vert
_{L^{2q}\left(  \Omega\right)  }.
\]
We choose $q=\infty$ and $q^{\prime}=1$ and use (\ref{Embedding}) to derive
$\left\Vert \kappa^{\prime}f\right\Vert \leq C\left\Vert \kappa^{\prime
}\right\Vert \left\Vert f\right\Vert _{H^{1}\left(  \Omega\right)  }$. This is
the upper bound with $C_{\kappa,1}=\left\Vert \kappa\right\Vert _{L^{\infty
}\left(  \Omega\right)  }+C\left\Vert \kappa^{\prime}\right\Vert $. To get the
lower bound we start with
\[
\rVert f\rVert_{H^{1}\left(  \Omega\right)  }=\rVert\frac{1}{\kappa}\kappa
f\rVert_{H^{1}\left(  \Omega\right)  }\leq\left\Vert \left(  \frac{1}{\kappa
}\right)  ^{\prime}\kappa f\right\Vert +\left\Vert \frac{1}{\kappa}\right\Vert
_{L^{\infty}\left(  \Omega\right)  }\left\Vert \left(  \kappa f\right)
^{\prime}\right\Vert \leq\left\Vert \left(  \frac{1}{\kappa}\right)  ^{\prime
}\kappa f\right\Vert +\frac{1}{\kappa_{\min}}\left\Vert \left(  \kappa
f\right)  ^{\prime}\right\Vert .
\]
We use the previous result to obtain%
\[
\left\Vert \left(  \frac{1}{\kappa}\right)  ^{\prime}\kappa f\right\Vert \leq
C\left\Vert \frac{\kappa^{\prime}}{\kappa^{2}}\right\Vert \left\Vert \kappa
f\right\Vert _{H^{1}\left(  \Omega\right)  }\leq\frac{C}{\kappa_{\min}^{2}%
}\left\Vert \kappa^{\prime}\right\Vert \left\Vert \kappa f\right\Vert
_{H^{1}\left(  \Omega\right)  }.
\]
Hence, the lower bound follows from%
\[
\rVert f\rVert_{H^{1}\left(  \Omega\right)  }\leq c_{\kappa,1}^{-1}\left\Vert
\kappa f\right\Vert _{H^{1}\left(  \Omega\right)  }\quad\text{with\quad
}c_{\kappa,1}:=\left(  1+C\frac{\left\Vert \kappa^{\prime}\right\Vert
_{L^{2}\left(  \Omega\right)  }}{\kappa_{\min}}\right)  ^{-1}\kappa_{\min}.
\]
Next, we consider the case $1/2\leq s<1$. For some $0<\varepsilon<1/2$, we
substitute in Lemma \ref{brezis}: $f\leftarrow\kappa$, $g\leftarrow f$,
$p\leftarrow2$, $q\leftarrow\varepsilon^{-1}$, $\theta\leftarrow1-\varepsilon
$, $t\leftarrow2-\frac{2\varepsilon}{2\varepsilon-1}$, $s\leftarrow\frac
{s}{1-\varepsilon}$ and obtain
\begin{equation}
\rVert\kappa f\rVert_{H^{s}\left(  \Omega\right)  }\leq C\left(  \rVert
\kappa\rVert_{L^{\infty}\left(  \Omega\right)  }\rVert f\rVert_{H^{s}\left(
\Omega\right)  }+\rVert f\rVert_{L^{1/\varepsilon}}\rVert\kappa\rVert
_{W^{\frac{s}{1-\varepsilon},2-\frac{2\varepsilon}{2\varepsilon-1}}\left(
\Omega\right)  }^{1-\varepsilon}\rVert\kappa\rVert_{L^{\infty}\left(
\Omega\right)  }^{\varepsilon}\right)  . \label{kappafslarge}%
\end{equation}
Sobolev's embedding theorem implies that there is a constant $C_{\varepsilon}
$ such that $\left\Vert f\right\Vert _{L^{1/\varepsilon}\left(  \Omega\right)
}\leq C_{\varepsilon}\left\Vert f\right\Vert _{H^{s}\left(  \Omega\right)  }$.
Hence, the upper estimate holds with%
\begin{equation}
C_{\kappa,s}:=C\left(  \rVert\kappa\rVert_{L^{\infty}\left(  \Omega\right)
}+C_{\varepsilon}\rVert\kappa\rVert_{W^{\frac{s}{1-\varepsilon},2-\frac
{2\varepsilon}{2\varepsilon-1}}\left(  \Omega\right)  }^{1-\varepsilon}%
\rVert\kappa\rVert_{L^{\infty}\left(  \Omega\right)  }^{\varepsilon}\right)
\quad\forall1/2\leq s<1. \label{Ckepsupper}%
\end{equation}
For the lower bound, we obtain%
\begin{equation}
\rVert f\rVert_{H^{s}\left(  \Omega\right)  }=\rVert\frac{1}{\kappa}\kappa
f\rVert_{H^{s}\left(  \Omega\right)  }\leq C_{\frac{1}{\kappa},s}\rVert\kappa
f\rVert_{H^{s}\left(  \Omega\right)  }\quad\text{so that\quad}c_{\kappa
,s}:=C_{\frac{1}{\kappa},s}^{-1}. \label{ckappas}%
\end{equation}
Next, we consider the case $0<s<\frac{1}{2}$ and observe%
\[
H^{s}(\Omega)\hookrightarrow L^{r}(\Omega),\quad\text{for}\quad1\leq
r\leq\dfrac{2}{1-2s}.
\]
Hence we have to restrict $\varepsilon$ in (\ref{kappafslarge}) to $\frac
{1}{2}-s\leq\varepsilon<1/2$. The constant $C_{\kappa,s}$ has the same form
(\ref{Ckepsupper}) while $\varepsilon$ therein must be chosen from the reduced
range. The same holds for the lower bound: it is of the same form
(\ref{ckappas}) while $\varepsilon$ therein must be chosen from the reduced range.

It remains to prove the estimate for the $H^{-s}\left(  \Omega\right)  $. For
$0<s\leq1$, we have
\begin{align*}
\rVert\kappa f\rVert_{H^{-s}\left(  \Omega\right)  }  &  =\sup_{\omega\in
H^{s}\left(  \Omega\right)  \backslash\left\{  0\right\}  }\dfrac{|(\kappa
f,\omega)_{L^{2}\left(  \Omega\right)  }|}{\rVert\omega\rVert_{{H}^{s}\left(
\Omega\right)  }}=\sup_{\omega\in H^{s}\left(  \Omega\right)  \backslash
\left\{  0\right\}  }\dfrac{|(f,\kappa\omega)_{L^{2}\left(  \Omega\right)  }%
|}{\rVert\omega\rVert_{{H}^{s}\left(  \Omega\right)  }}\\
&  \leq\left\Vert f\right\Vert _{H^{-s}\left(  \Omega\right)  }\sup_{\omega\in
H^{s}\left(  \Omega\right)  \backslash\left\{  0\right\}  }\dfrac{\left\Vert
\kappa\omega\right\Vert _{H^{s}\left(  \Omega\right)  }}{\rVert\omega
\rVert_{{H}^{s}\left(  \Omega\right)  }}\leq C_{\kappa,s}\left\Vert
f\right\Vert _{H^{-s}\left(  \Omega\right)  }.
\end{align*}
To get the lower bound, we notice that
\begin{align*}
\rVert f\rVert_{H^{-s}\left(  \Omega\right)  }  &  =\rVert\frac{1}{\kappa
}\kappa f\rVert_{H^{-s}\left(  \Omega\right)  }=\sup_{\omega\in H^{s}\left(
\Omega\right)  \backslash\left\{  0\right\}  }\dfrac{|(\kappa f,\frac
{1}{\kappa}\omega)_{L^{2}\left(  \Omega\right)  }|}{\rVert\omega\rVert
_{H^{s}\left(  \Omega\right)  }}\\
&  \leq\rVert\kappa f\rVert_{H^{-s}\left(  \Omega\right)  }\sup_{\omega\in
H^{s}\left(  \Omega\right)  \backslash\left\{  0\right\}  }\dfrac{\rVert
\frac{1}{\kappa}\omega\rVert_{H^{s}\left(  \Omega\right)  }}{\rVert
\omega\rVert_{H^{s}\left(  \Omega\right)  }}\leq C_{\frac{1}{\kappa},s}%
\rVert\kappa f\rVert_{H^{-s}\left(  \Omega\right)  }.
\end{align*}%
\endproof

This norm equivalence allows us to generalize the class of kernel functions in
Abel's integral operator.

\begin{definition}
Let $s\in\left[  -1,1\right]  $ and $\Omega=\left(  0,1\right)  $. The class
of $s$-admissible kernel functions $\mathcal{A}\left(  s\right)  $ consists of
functions $K:\Omega\times\Omega\rightarrow\mathbb{R}$ such that there exists a
sequence $\left(  \psi_{n}\right)  _{n}$ of functions in the multiplier class
$\mathcal{M}\left(  s\right)  $ such that the following conditions hold:

\begin{enumerate}
\item The kernel of integral equation has the representation
\begin{equation}
K(x,y)=\sum_{n=1}^{\infty}\sum_{m=1}^{\infty}d_{n,m}\psi_{n}(x)\psi
_{m}(y)\qquad\forall x,y\in\Omega\quad\text{a.e.} \label{repK1}%
\end{equation}
for some coefficients $\left(  d_{n,m}\right)  _{n,m=1}^{\infty}$.

\item The constants $C_{n,s}$, $c_{n,s}$ in the multiplier estimates
\begin{equation}
c_{n,s}\rVert f\rVert_{H^{s}\left(  \Omega\right)  }\leq\rVert\psi_{n}%
f\rVert_{H^{s}\left(  \Omega\right)  }\leq C_{n,s}\rVert f\rVert_{H^{s}\left(
\Omega\right)  }\quad\forall f\in H^{s}\left(  \Omega\right)  \label{cnsCns}%
\end{equation}
satisfy:

\begin{enumerate}
\item there is a positive constant $C_{s}<\infty$ such that%
\[
\sum_{n,m=1}^{\infty}\vert d_{n,m}\vert C_{n,s}C_{m,s}\leq C_{s}^{2},
\]

\item there is a constant $\tilde{\gamma}$ such that
\begin{equation}
\gamma\sum_{n=1}^{\infty}d_{n,n}c_{n,s}^{2}-C_{c}\sum_{\substack{n,m=1\\n\neq
m}}^{\infty}|d_{n,m}|C_{n,s}C_{m,s}\geq\tilde{\gamma}>0,
\label{defgammatildeconst}%
\end{equation}
where $\gamma$ is as in (\ref{ellineq}) and $C_{c}$ as in
(\ref{continuityestimate}).
\end{enumerate}
\end{enumerate}
\end{definition}

\begin{theorem}
\label{main1}Suppose that the kernel function $K$ in Abel's integral equation
is $s$-admissible for $s=-\alpha/2$. Then $a_{K,\alpha}(f,g)$ is a continuous
and $H^{-\alpha/2}\left(  \Omega\right)  $-elliptic sesquilinear form.
\end{theorem}

%

\proof
Let $s=-\alpha/2$. We have%
\begin{align*}
\rvert a_{K,\alpha}(f,g)\rvert &  \leq\sum_{n,m=1}^{\infty}\rvert
d_{n,m}a_{1,\alpha}(\psi_{m}f,\psi_{n}g)\rvert\\
&  \leq C_{c}\sum_{n,m=1}^{\infty}|d_{n,m}|\rVert\psi_{m}f\rVert_{H^{s}\left(
\Omega\right)  }\rVert\psi_{n}g\rVert_{H^{s}\left(  \Omega\right)  }\\
&  \leq C_{c}\left(  \sum_{n,m=1}^{\infty}|d_{n,m}|C_{n,s}C_{m,s}\right)
\rVert f\rVert_{H^{s}\left(  \Omega\right)  }\rVert g\rVert_{H^{s}\left(
\Omega\right)  }\\
&  \leq\tilde{C}_{c}\rVert f\rVert_{H^{s}\left(  \Omega\right)  }\rVert
g\rVert_{H^{s}\left(  \Omega\right)  }\quad\text{with\quad}\tilde{C}%
_{c}:=C_{c}C_{s}^{2}.
\end{align*}
To obtain the $H^{s}\left(  \Omega\right)  $-ellipticity, we observe that%
\begin{align*}
a_{K,\alpha}(f,f)  &  =\sum_{n,m=1}^{\infty}d_{n,m}a_{1,\alpha}(\psi_{n}%
f,\psi_{m}f)\\
&  =\sum_{n=1}^{\infty}d_{n,n}a_{1,\alpha}(\psi_{n}f,\psi_{n}f)+\sum
_{\substack{n,m=1\\n\neq m}}^{\infty}d_{n,m}a_{1,\alpha}(\psi_{n}f,\psi_{m}f).
\end{align*}
We employ the coercivity and continuity of $a_{1,\alpha}$ to derive%
\[
\operatorname{Re}a_{K,\alpha}(f,f)\geq\gamma\sum_{n=1}^{\infty}d_{n,n}%
\rVert\psi_{n}f\rVert_{H^{s}\left(  \Omega\right)  }^{2}-C_{c}\sum
_{\substack{n,m=1\\n\neq m}}^{\infty}|d_{n,m}|\left\vert a_{1,\alpha}(\psi
_{n}f,\psi_{m}f)\right\vert .
\]
The estimates for the multipliers $\psi_{n}$ lead to%
\begin{equation}
\operatorname{Re}a_{K,\alpha}(f,f)\geq\left(  \gamma\sum_{n=1}^{\infty}%
d_{n,n}c_{n,s}^{2}-C_{c}\sum_{\substack{n,m=1\\n\neq m}}^{\infty}%
|d_{n,m}|C_{n,s}C_{m,s}\right)  \rVert f\rVert_{H^{s}\left(  \Omega\right)
}^{2}\geq\tilde{\gamma}\rVert f\rVert_{H^{s}\left(  \Omega\right)  }^{2}.
\label{ellcoeraKalpha}%
\end{equation}
The assumptions on the summability of the constants $c_{n,s}$, $C_{n,s}$ lead
to the assertion.%
\endproof

\section{Discretization of Abel-type Integral Equations \label{Sec3}}

In order to solve the Abel-type integral equation (\ref{gal1}), (\ref{var1})
numerically we discretize the continuous problem (\ref{var1}) by a Galerkin
finite element method. For this we introduce the piecewise polynomial finite
element spaces. Let a set of mesh points $\mathcal{N}=\left(  x_{i}\right)
_{i=0}^{N}$ be given%
\[
0=x_{0}<x_{1}<\ldots<x_{N}=1
\]
which induces a mesh $\mathcal{T}=\left\{  \tau_{i}:1\leq i\leq N\right\}  $
on $\Omega$, where $\tau_{i}=\left[  x_{i-1},x_{i}\right]  $. The length of a
subinterval $\tau\in\mathcal{T}$ is denoted by $h_{\tau}$ and the maximal mesh
width by $h:=\max\left\{  h_{\tau}:\tau\in\mathcal{T}\right\}  $. The
variation of the lengths of neighboring intervals is controlled by the
constant%
\[
C_{\mathcal{T}}:=\max\left\{  \frac{h_{\tau}}{h_{\sigma}}:\forall\tau
,\sigma\in\mathcal{T}\text{ with }\tau\cap\sigma\neq\emptyset\right\}  .
\]
The piecewise polynomial function space of degree $m\in\mathbb{N}_{0}$ on
$[0,1]$ is given by
\begin{equation}
S_{\mathcal{T}}^{m}:=\left\{
\begin{array}
[c]{ll}%
\{v\in L^{\infty}(\Omega):\left.  v\right\vert _{\tau}\in\mathbb{P}_{0}\left(
\tau\right)  ,\forall\tau\in\mathcal{T}\} & m=0,\\
\{v\in C(\Omega):\left.  v\right\vert _{\tau}\in\mathbb{P}_{m}\left(
\tau\right)  ,\forall\tau\in\mathcal{T}\} & m\geq1.
\end{array}
\right.
\end{equation}
Here $\mathbb{P}_{m}\left(  \tau\right)  $ denote the space of all univariate
polynomials on $\tau$ of maximal degree $m$. The nodal points are given by%
\[
\mathcal{N}_{m}:=\left\{
\begin{array}
[c]{ll}%
\left\{  \frac{x_{i}+x_{i-1}}{2}:1\leq i\leq N\right\}  & m=0,\\
\left\{  \xi_{i,j}:=x_{i-1}+j\frac{x_{i}-x_{i-1}}{m}\quad1\leq i\leq N\text{,
}0\leq j\leq m-1\right\}  \cup\left\{  1\right\}  & m\geq1,
\end{array}
\right.
\]
so that the dimension of $S_{\mathcal{T}}^{m}$ is $M:=N$ for $m=0$ and
$M:=Nm+1$ for $m\geq1$. We choose the usual Lagrange basis functions
$b_{i}^{(m)}$ of $S_{\mathcal{T}}^{m}$ and write $b_{i}$ short for
$b_{i}^{(m)}$ if the polynomial degree is clear from the context. The Galerkin
method is given by replacing the infinite-dimensional space $H^{-\alpha
/2}\left(  \Omega\right)  $ in (\ref{var1}) by the finite dimensional subspace
$S:=S_{\mathcal{T}}^{m}$:
\begin{equation}
\text{Find\: }f_{S}\in S\ \text{such that\qquad}a_{K,\alpha}(f_{S}%
,\varphi)=G(\varphi)\quad\forall\varphi\in S. \label{discprob}%
\end{equation}
For the computation of $f_{S}$ one introduces the representation of
(\ref{discprob}) with respect to the basis $(b_{i})_{i=0}^{M}$ of $S$. Let
\[
a_{ij}=(A_{K,\alpha}b_{j},b_{i})=\Gamma(\alpha)^{-1}\int_{0}^{1}{b_{i}%
(x)\int_{0}^{x}{(x-y)^{\alpha-1}K(x,y)b_{j}(y)dy}dx},
\]%
\begin{equation}
r_{i}=(g,b_{i})=\int_{0}^{1}{b_{i}g}, \label{rhs}%
\end{equation}
with the \textit{system matrix} $\mathbf{A}=(a_{ij})_{i,j=1}^{M}\in
\mathbb{R}^{M\times M}$ and the \textit{right-hand side vector} $\mathbf{r}%
=(r_{i})_{i=1}^{M}\in\mathbb{R}^{M}$. Then, the basis representation of
(\ref{discprob}) is: Find $\mathbf{f}_{S}\in\mathbb{R}^{M}$ such that
\begin{equation}
\mathbf{Af}_{S}\mathbf{=r} \label{lse}%
\end{equation}
and the solution of (\ref{discprob}) is given by $f_{S}=\sum_{i=1}^{M}%
f_{S,i}b_{i}$.

\section{Convergence Analysis \label{Sec4}}

\subsection{Discretization Error}

In order to estimate the discretization error of the Galerkin discretization
we employ C\'{e}a's lemma.

\begin{theorem}
[C\'{e}a]In Abel's integral operator $A_{K,\alpha}$, let $\alpha\in\left(
0,1\right)  $ and $K\in\mathcal{A}\left(  -\frac{\alpha}{2}\right)  $. For
some $G\in H^{\alpha/2}\left(  \Omega\right)  $, let $f\in H^{-\alpha
/2}\left(  \Omega\right)  $ be the exact solution of (\ref{var1}). Then also
the Galerkin discretization (\ref{discprob}) has a unique solution $f_{S}\in
S$ which satisfies the quasi-optimal error estimate%
\begin{equation}
\Vert f-f_{S}\Vert_{H^{-\alpha/2}\left(  \Omega\right)  }\leqslant\frac
{\tilde{C}_{c}}{\tilde{\gamma}}\inf\limits_{v\in S}\Vert f-v\Vert
_{H^{-\alpha/2}\left(  \Omega\right)  }, \label{discerror1}%
\end{equation}
where $\tilde{\gamma},\tilde{C}_{c}$ are the ellipticity and continuity
constants of the form $a_{K,\alpha}$ as in the proof of Theorem \ref{main1}.
\end{theorem}

To derive convergence \textit{rates} we investigate the error term $\left.
\inf_{v\in S_{\mathcal{T}}^{m}}\left\Vert f-v\right\Vert _{H^{-\alpha
/2}\left(  \Omega\right)  }\right.  $. It is well known (see, e.g.,
{\cite{cars, sauter11}) }that for sufficiently smooth solution $f\in
H^{m+1}\left(  \Omega\right)  $ it holds%
\[
\Vert f-f_{S}\Vert_{H^{-\alpha/2}\left(  \Omega\right)  }\leqslant\frac
{\tilde{C}_{c}}{{\tilde{\gamma}}}\min\limits_{v\in S_{\mathcal{T}}^{m}}\Vert
f-v\Vert_{H^{-\alpha/2}\left(  \Omega\right)  }\leqslant\frac{\tilde{C}_{c}%
}{{\tilde{\gamma}}}Ch^{m+1+\alpha/2}\Vert f\Vert_{H^{m+1}\left(
\Omega\right)  }.
\]

\subsection{Perturbation (Quadrature) \label{SubSecPert}}

To derive a fully discrete method we apply numerical quadrature to approximate
the integrals in the system matrix. We develop the quadrature error analysis
for analytic $K$, more precisely, we assume that there exists constants
$C_{K}$ and $\Lambda_{K}$ such that%
\begin{equation}
\left\Vert K\right\Vert _{C^{n}\left(  \Omega\times\Omega\right)  }\leq
C_{K}\Lambda_{K}^{n}n!,\qquad\forall n\in\mathbb{N}_{0}. \label{analyticK}%
\end{equation}
To reduce technicalities we assume that $\Lambda_{K}\geq2$. We introduce the
function $w_{i,j}(x,y):=K(x,y)(x-y)^{\alpha-1}b_{i}(x)b_{j}(y)$ and
reformulate the integral for a matrix entry $a_{ij}=I\left(  w_{i,j}\right)
$:
\begin{align}
\label{integrals}I\left(  w_{i,j}\right)   &  =\frac{1}{\Gamma(\alpha)}%
\int_{0}^{1}{\int_{0}^{x}w_{i,j}(x,y){dy}dx}\\
&  =\frac{1}{\Gamma(\alpha)}\sum_{\tau\subset\operatorname*{supp}(b_{i})}%
{\sum_{\sigma\subset\operatorname*{supp}(b_{j})}{\int_{\tau}{\int_{\sigma
\cap\lbrack0,x]}w_{i,j}(x,y){dy}dx}}}\nonumber
\end{align}
with $\alpha\in(0,1)$. Since all the above functions are real we may omit the
complex conjugation on the second argument in the $L^{2}$ scalar product. We
compute the integral in (\ref{integrals}) in different ways depending on the
relative location of $\tau$ compared to $\sigma$:

\begin{itemize}
\item if $\tau=\sigma$ we employ simplex coordinates (cf. Section \ref{case2}),

\item if $\tau$ lies to the left of $\sigma$ it is easy to see that integral
value is 0,

\item in the case that $\tau$ lies to the right of $\sigma$ we apply
tensor-Gauss quadrature (cf. Section \ref{case1}).
\end{itemize}

\subsubsection{The Case $\tau\neq\sigma$}

\label{case1} The $n\times n$-tensor-Gauss quadrature for a function $f$ on an
interval $\tau\times\sigma$ is given by
\[
\left(  Q_{\tau}^{n}\otimes Q_{\sigma}^{n}\right)  (f)=\sum_{k=0}^{n-1}%
{\sum_{\ell=0}^{n-1}{\omega_{k}^{\tau,G}\omega_{\ell}^{\sigma,G}f(\xi
_{k}^{\tau,G},\xi_{\ell}^{\tau,G})}},
\]
with $\omega_{k}^{\tau,G},\xi_{k}^{\tau,G}$ denoting the weights and abscissae
for $Q_{\tau}^{n}$, the $n$-point Gauss-Legendre quadrature method scaled to
$\tau$. The exact integral value is given by%
\[
\left(  I_{\tau}\otimes I_{\sigma}\right)  (w_{i,j})=\int_{\tau}{\int_{\sigma
}{(x-y)^{\alpha-1}K(x,y)b_{i}(x)b_{j}(y)dy}dx}.
\]
We have the following error estimate \cite[Chap. 5]{sauter11}:%
\begin{align}
\label{tenserr}|E^{n}f|  &  =|(I_{\tau}\otimes I_{\sigma}-Q_{\tau}^{n}\otimes
Q_{\sigma}^{n})(f)|\\
&  \leqslant\max_{t\in\tau}|(I_{\sigma}-Q_{\sigma}^{n})f(t,\cdot)|+\max
_{t\in\sigma}|(I_{\tau}-Q_{\tau}^{n})f(\cdot,t)|\nonumber\\
&  =\max_{t\in\tau}|E_{\sigma}^{n}f(t,\cdot)|+\max_{t\in\sigma}|E_{\tau}%
^{n}f(\cdot,t)|,\nonumber
\end{align}
where $I_{\tau}$ and $I_{\sigma}$ are the integrals with respect to the first
and second variable, respectively. We first investigate the error term
$\max_{t\in\tau_{i}}|E_{\sigma}^{n}f(t,\cdot)|$ by the following lemma.

\begin{lemma}
\label{LemQuad1}Let $b_{j}$, $j\in\{1,\ldots,M\}$, be the Lagrange basis of
$S_{\mathcal{T}}^{m}$ and let $K$ satisfy (\ref{analyticK}). Assume that the
number of quadrature points satisfy $2n>m$. Let $x\in\Omega$, $\sigma
\subset\operatorname*{supp}b_{j}$ satisfy $\operatorname*{dist}\left(
x,\sigma\right)  >0$.
We have%
\begin{equation}
\left\vert (I_{\sigma}-Q_{\sigma}^{n})\left(  \left(  x-\cdot\right)
^{\alpha-1}K\left(  x,\cdot\right)  b_{j}\right)  \right\vert \leq\frac
{2C_{K}\operatorname*{e}\Lambda_{K}}{\operatorname*{dist}\nolimits^{1-\alpha
}\left(  x,\sigma\right)  }\left(  \frac{\Lambda_{K}}{2}\frac{h_{\sigma}%
}{\operatorname*{dist}\left(  x,\sigma\right)  }\right)  ^{2n}\left\Vert
b_{j}\right\Vert _{C^{m}\left(  \sigma\right)  }. \label{gausserrcol}%
\end{equation}

\end{lemma}%

\proof
Since the exactness degree of an $n$-point Gaussian quadrature method is
$2n-1$, it holds (cf. \cite{gautschi12})%
\begin{equation}
\left\Vert \left(  I_{\sigma}-Q_{\sigma}^{n}\right)  \left(  (x-\cdot
)^{\alpha-1}K(x,\cdot)b_{j}\right)  \right\Vert _{L^{\infty}\left(
\sigma\right)  }\leqslant2\inf_{p\in\mathbb{P}_{2n-1}}\left\Vert \left(
{x-\cdot}\right)  {^{\alpha-1}K(x,\cdot)b_{j}-p}\right\Vert {_{L^{\infty
}\left(  \sigma\right)  }}. \label{taylor}%
\end{equation}
We compute the $2n^{th}$ derivative $\partial^{2n}\left(  (x-\cdot)^{\alpha
-1}K(x,\cdot)b_{j}\right)  $ to estimate the remainder in a Taylor expansion%
\[
\left\vert \partial^{2n}\left(  (x-\cdot)^{\alpha-1}K(x,\cdot)b_{j}\right)
\right\vert =\left\vert \sum_{k=0}^{2n}\binom{2n}{k}\left(  \partial^{k}%
b_{j}\right)  \partial^{2n-k}\left(  (x-\cdot)^{\alpha-1}K(x,\cdot)\right)
\right\vert .
\]
For $\mu\in\{0,\ldots,2n\}$, the derivative in the last term can be estimated
by%
\begin{align*}
\left\Vert \partial^{\mu}\left(  (x-\cdot)^{\alpha-1}K(x,\cdot)\right)
\right\Vert _{L^{\infty}\left(  \sigma\right)  }  &  =\left\Vert \sum_{\nu
=0}^{\mu}{\binom{\mu}{\nu}\partial^{\nu}K(x,\cdot)\partial^{\mu-\nu}%
(x-\cdot)^{\alpha-1}}\right\Vert _{L^{\infty}\left(  \sigma\right)  }\\
&  \leqslant\left\Vert \sum_{\nu=0}^{\mu}{\binom{\mu}{\nu}\partial^{\nu
}K(x,\cdot)(\mu-\nu)!(x-\cdot)^{\alpha-1-(\mu-\nu)}}\right\Vert _{L^{\infty
}\left(  \sigma\right)  }\\
&  \leqslant\left\Vert \mu!\left(  x-\cdot\right)  ^{\alpha-1-\mu}\sum_{\nu
=0}^{\mu}\frac{{\left\vert \partial^{\nu}K(x,\cdot)\right\vert }}{\nu
!}\right\Vert _{L^{\infty}\left(  \sigma\right)  },
\end{align*}
where we used $\operatorname*{dist}^{\nu}\left(  x,\sigma\right)  \leq1$.
Since
\[
\sum_{\nu=0}^{\mu}\frac{\left\Vert \partial^{\nu}K(x,\cdot)\right\Vert
_{L^{\infty}\left(  \sigma\right)  }}{\nu!}\leq\sum_{\nu=0}^{\mu}\frac
{C_{K}\nu!\Lambda_{K}^{\nu}}{\nu!}\leq C_{K}\Lambda_{K}^{\mu+1},
\]
and $\sum_{k=0}^{\infty}\frac{1}{k!}=\operatorname*{e}$, we obtain
\begin{align*}
\left\Vert \partial^{2n}\left(  (x-\cdot)^{\alpha-1}K(x,\cdot)b_{j}\right)
\right\Vert _{L^{\infty}\left(  \sigma\right)  }  &  \leqslant\sum_{k=0}%
^{2n}{\frac{(2n)!}{k!}}\left\Vert {\left\vert \partial^{k}b_{j}\right\vert
(x-\cdot)^{\alpha-1-2n}\sum_{\nu=0}^{2n-k}\frac{{\left\vert \partial^{\nu
}K(x,\cdot)\right\vert }}{\nu!}}\right\Vert _{L^{\infty}\left(  \sigma\right)
}\\
&  \leqslant C_{K}\operatorname*{e}\Lambda_{K}^{2n+1}(2n)!\Vert b_{j}%
\Vert_{C^{m}(\sigma)}\operatorname*{dist}\left(  x,\sigma\right)
^{\alpha-1-2n}.
\end{align*}
This leads to%
\begin{align*}
&  \left\Vert \left(  I_{\sigma}-Q_{\sigma}^{n}\right)  \left(  (x-\cdot
)^{\alpha-1}K(x,\cdot)b_{j}\right)  \right\Vert _{L^{\infty}\left(
\sigma\right)  }\\
&  \qquad\qquad\leqslant2\frac{\left(  h_{\sigma}/2\right)  ^{2n}}{\left(
2n\right)  !}\sum_{k=0}^{2n}\frac{\left(  2n\right)  !}{k!}\left\Vert \left(
\partial^{k}b_{j}\right)  \left(  x-\cdot\right)  ^{\alpha-1-\left(
2n-k\right)  }\sum_{\nu=0}^{2n-k}\frac{{\left\vert \partial^{\nu}%
K(x,\cdot)\right\vert }}{\nu!}\right\Vert _{L^{\infty}\left(  \sigma\right)
}\\
&  \qquad\qquad\leq2\left(  \frac{h_{\sigma}}{2}\right)  ^{2n}\left\Vert
b_{j}\right\Vert _{C^{m}\left(  \sigma\right)  }C_{K}\Lambda_{K}%
^{2n+1}\operatorname*{dist}\nolimits^{\alpha-1-2n}\left(  x,\sigma\right)
\sum_{k=0}^{2n}\frac{1}{k!}\\
&  \qquad\qquad\leq\frac{2C_{K}\operatorname*{e}\Lambda_{K}}%
{\operatorname*{dist}\nolimits^{1-\alpha}\left(  x,\sigma\right)  }\left(
\frac{h_{\sigma}\Lambda_{K}}{2\operatorname*{dist}\left(  x,\sigma\right)
}\right)  ^{2n}\left\Vert b_{j}\right\Vert _{C^{m}\left(  \sigma\right)  }%
\end{align*}
and we arrive at (\ref{gausserrcol}).
\endproof

We now investigate the error term $\max_{t\in\sigma}|E_{\tau}^{n}f(\cdot,t)|$
by the following lemma.

\begin{lemma}
\label{LemQuad1.5}Let $b_{i}$, $i\in\{1,\ldots, M\}$, be the Lagrange basis
for $S_{\mathcal{T}}^{m}$ and let $K$ satisfy (\ref{analyticK}) and $2n>$
$m$.
Then we have%
\begin{equation}
\left\vert (I_{\tau}-Q_{\tau}^{n})\left(  (\cdot-y)^{\alpha-1}K(\cdot
,y)b_{i}\right)  \right\vert \leq\frac{2C_{K}\operatorname*{e}\Lambda_{K}%
}{\operatorname*{dist}\nolimits^{1-\alpha}\left(  x, \tau\right)  }\left(
\frac{\Lambda_{K}}{2}\frac{h_{\tau}}{\operatorname*{dist}\left(  y,
\tau\right)  }\right)  ^{2n}\left\Vert b_{i}\right\Vert _{C^{m}\left(
\tau\right)  }. \label{erri}%
\end{equation}

\end{lemma}

The proof is a repetition of the arguments used in the proof of Lemma
\ref{LemQuad1} and skipped here. The proof of the following theorem follows
from Lemma \ref{LemQuad1} and \ref{LemQuad1.5} via a straightforward tensor
argument (cf. (\ref{tenserr})).

\begin{theorem}
\label{TheoLocQuad}Let $b_{i},b_{j}$, $i,j\in\{1,\ldots, M\}$, be Lagrange
basis functions of the piecewise polynomial space $S_{\mathcal{T}}^{m}$ and
let $K$ satisfy (\ref{analyticK}). Assume $2n>m$. For $\tau\subset
\operatorname*{supp}(b_{i}),\sigma\subset\operatorname*{supp}(b_{j})$ with
$\operatorname*{dist}\left(  \tau,\sigma\right)  >0$ it holds%
\begin{equation}
\left\vert (I_{\tau}\otimes I_{\sigma}-Q_{\tau}^{n}\otimes Q_{\sigma}%
^{n})\left(  w_{i,j}\right)  \right\vert \leq\frac{2C_{K}\operatorname*{e}%
\Lambda_{K}}{\operatorname*{dist}\nolimits^{1-\alpha}\left(  \tau
,\sigma\right)  }\left(  \frac{\Lambda_{K}}{2}\frac{\max\left\{  h_{\tau
},h_{\sigma}\right\}  }{\operatorname*{dist}\left(  \tau,\sigma\right)
}\right)  ^{2n}\left(  \left\Vert b_{i}\right\Vert _{C^{m}\left(  \tau\right)
}+\left\Vert b_{j}\right\Vert _{C^{m}\left(  \sigma\right)  }\right)  .
\label{gauserrgal}%
\end{equation}

\end{theorem}

\begin{remark}
\label{RemExact}The upper bound for the quadrature error in Theorem
\ref{TheoLocQuad} grows if $\max\left\{  h_{\tau},h_{\sigma}\right\}
>\frac{2}{\Lambda_{K}}\operatorname*{dist}(\tau,\sigma)$. This is an artifact
of the proof. If one employs the theory of derivative-free error quadrature
error estimates via complex analysis of analytic integrands (cf. \cite{DR},
\cite{Tref_Approx}) one obtains exponential convergence as long as
$\operatorname*{dist}\left(  \tau,\sigma\right)  >0$. Here, we employed the
simpler classical theory in order to reduce technicalities. For the numerical
experiments we have used the tensor Gauss-Legendre formulae for all pairs of
intervals $\tau,\sigma$ which have positive distance.
\end{remark}

\subsubsection{The Case $\tau=\sigma$}

\label{case2} In order to approximate the integral
\[
\int_{\tau}{\int_{\tau\cap\lbrack0,x]}{(x-y)^{\alpha-1}K(x,y)b_{i}%
(x)b_{j}(y)dy}dx}%
\]
we first transform it to the unit cube via simplex coordinates $\left(
\xi,\eta\right)  \rightarrow\left(  \xi,\xi\eta\right)  $. The determinant of
the Jacobean equals $\xi$. Let $\tau=[a,a+h_{\tau}]$; then%
\begin{align*}
&  \frac{1}{\Gamma(\alpha)}\int_{a}^{a+h_{\tau}}{\int_{a}^{x}{(x-y)^{\alpha
-1}K(x,y)b_{i}(x)b_{j}(y)dy}dx}\\
&  \qquad\qquad=\frac{1}{\Gamma(\alpha)}\int_{0}^{h_{\tau}}{\int_{0}%
^{x}{(x-y)^{\alpha-1}K(x+a,y+a)b_{i}(x+a)b_{j}(y+a)dy}dx}\\
&  \qquad\qquad=\frac{1}{\Gamma(\alpha)}\int_{0}^{h_{\tau}}{\int_{0}^{1}%
{\xi^{\alpha}(1-\eta)^{\alpha-1}K(\xi+a,\xi\eta+a)b_{i}(\xi+a)b_{j}(\xi
\eta+a)d\eta}d\xi}\\
&  \qquad\qquad=\frac{h_{\tau}^{\alpha+1}}{\Gamma(\alpha)}\int_{0}^{1}%
{\int_{0}^{1}{\xi^{\alpha}(1-\eta)^{\alpha-1}K(h_{\tau}\xi+a,h_{\tau}\xi
\eta+a)b_{i}(h_{\tau}\xi+a)b_{j}(h_{\tau}\xi\eta+a)d\eta}d\xi.}%
\end{align*}
For the $\xi$-integration we employ Gauss-Jacobi quadrature with weight
$\xi^{\alpha}$ and for the $\eta$-integration we employ Gauss-Jacobi
quadrature with weight $(1-\eta)^{\alpha-1}$. For an error analysis for these
method we refer to \cite[Section 5.3.2]{sauter11}.

\subsubsection{Right-Hand Side}

The right-hand side of the Galerkin method can be approximated with an
$n$-point Gaussian quadrature method:
\begin{equation}
r_{i}=\sum_{\tau\subset\operatorname*{supp}(b_{i})}{\int_{\tau}{gb_{i}}%
}\approx\sum_{\tau\subset\operatorname*{supp}(b_{i})}{\sum_{k=0}^{n-1}%
{\omega_{k}^{\tau,G}g(\xi_{k}^{\tau,G})b_{i}(\xi_{k}^{\tau,G})}}
\label{rhserr0}%
\end{equation}
The error of the right-hand side is estimated in the following lemma. We
assume that the right-hand side $g$ is analytic, more precisely, there exist
constants $C_{g}$ and $\Lambda_{g}\geq2$ such that%
\begin{equation}
\left\Vert g\right\Vert _{C^{n}\left(  \Omega\right)  }\leq C_{g}\Lambda
_{g}^{n}n!\qquad\forall n\in\mathbb{N}_{0}. \label{ganalytic}%
\end{equation}

\begin{lemma}
Let $b_{i}$, $i\in\{1,\ldots,M\}$, denote the Lagrange basis of
$S_{\mathcal{T}}^{m}$. Assume $2n>m$ and $g$ satisfies (\ref{ganalytic}). Let
$\tau\subset\operatorname*{supp}(b_{i})$. The error of the $n$-point Gauss
quadrature method is given by
\begin{equation}
\left\vert (I_{\tau}-Q_{\tau}^{n})(b_{i}g)\right\vert \leqslant2\Lambda
_{g}\left(  \frac{\Lambda_{g}}{2}h_{\tau}\right)  ^{2n}\Vert b_{i}\Vert
_{C^{m}(\tau)}. \label{gausserrrhs}%
\end{equation}

\end{lemma}%

\proof
We compute the $2n^{th}$ derivative of the integrand as in the previous
section:
\begin{align*}
\left\Vert \partial^{2n}\left(  b_{i}g\right)  \right\Vert _{L^{\infty}\left(
\tau\right)  }  &  =\left\Vert \sum_{k=0}^{2n}\binom{2n}{k}\left(
\partial^{k}b_{i}\right)  \left(  \partial^{2n-k}g\right)  \right\Vert
_{L^{\infty}\left(  \tau\right)  }\\
&  =\left\Vert b_{i}g^{(2n)}+\ldots+\binom{2n}{m}b_{i}^{\left(  m\right)
}g^{\left(  2n-m\right)  }\right\Vert _{L^{\infty}\left(  \tau\right)  }\\
&  \leqslant(2n)!\Vert b_{i}\Vert_{C^{m}(\tau)}\sum_{k=0}^{2n}\frac{\Vert
g\Vert_{C^{k}(\tau)}}{k!}\\
&  \leq C_{g}(2n)!\Vert b_{i}\Vert_{C^{m}(\tau)}\Lambda_{g}^{2n+1}%
\end{align*}
This allows us to estimate the error via a Taylor argument%
\[
|(I_{\tau}-Q_{\tau}^{n})(b_{i}g)|\leqslant2C_{g}(2n)!\Vert b_{i}\Vert
_{C^{m}(\tau)}\Lambda_{g}^{2n+1}\left(  \frac{(h_{\tau}/2)^{2n}}%
{(2n)!}\right)  =2\Lambda_{g}\left(  \frac{\Lambda_{g}}{2}h_{\tau}\right)
^{2n}\Vert b_{i}\Vert_{C^{m}(\tau)}.
\]
\endproof

\subsection{The Influence of the Quadrature on the Discretization
Error\label{SecInfluence}}

The approximation of the entries $\tilde{a}_{i,j}$, $\tilde{r}_{i}$ of the
system matrix and of the right-hand side vector by numerical quadrature leads
to a \textquotedblleft perturbed" linear system of equations
\begin{equation}
\tilde{\mathbf{A}}\tilde{\mathbf{f}}_{S}=\tilde{\mathbf{r}}, \label{persoe}%
\end{equation}
which can be translated to the following variational problem: Find $\tilde
{f}_{S}\in S$ such that%
\begin{equation}
\tilde{a}_{K,\alpha}(\tilde{f}_{S},v)=\tilde{G}(v)\quad\forall v\in S.
\label{pervarprob}%
\end{equation}
Here, $\tilde{a}_{K,\alpha}:S\times S\rightarrow\mathbb{R}$ and $\tilde
{G}:S\rightarrow\mathbb{R}$ are defined for $u=\sum_{i=1}^{M}\alpha_{i}b_{i}$
and $v=\sum_{i=1}^{M}\beta_{i}b_{i}$ by%
\[
\tilde{a}_{K,\alpha}\left(  u,v\right)  =\sum_{i,j=1}^{M}\beta_{i}\tilde
{a}_{i,j}\alpha_{j}\quad\text{and\quad}\tilde{G}\left(  v\right)  =\sum
_{i=1}^{M}\tilde{r}_{i}\beta_{i}.
\]

We will analyze the error $\Vert f-\tilde{f}_{S}\Vert_{H^{-\alpha/2}\left(
\Omega\right)  }$ and the influence of the quadrature on the discretization
error. To reduce technicalities we assume that the mesh is quasi-uniform,
i.e., there exists a constant $C_{\operatorname*{qu}}$ such that%
\[
C_{\operatorname*{qu}}:=\max\left\{  \frac{h}{h_{\tau}}:\tau\in\mathcal{T}%
\right\}  .
\]
First we rewrite the sesquilinear form $a(u,v)$ for $u,v\in S_{\mathcal{T}%
}^{m},m\in\mathbb{N}_{0}$
\begin{align}
a_{K,\alpha}(u,v)  &  =\frac{1}{\Gamma(\alpha)}\int_{0}^{1}{u(x)\int_{0}%
^{x}{(x-y)^{\alpha-1}K(x,y)v(y)dy}dx}\label{auf}\\
&  =\frac{1}{\Gamma(\alpha)}\sum_{i=1}^{M}\sum_{j=1}^{i}u_{i}{{v_{j}}%
\sum_{\tau\subset\operatorname*{supp}(b_{i})}\sum_{\sigma\subset
\operatorname*{supp}(b_{j})}{\int_{\tau}{{{{\int_{\sigma\cap\lbrack0,x]}}}%
}b_{i}(x){{{\frac{K(x,y)b_{j}(y)}{(x-y)^{1-\alpha}}dy}}}dx.}}}\nonumber
\end{align}
This leads to the definition $I_{\tau\times\sigma}^{i,j}:=\int_{\tau}%
{\int_{\sigma\cap\left[  0,x\right]  }b_{i}(x){b_{j}(y)K(x,y)(x-y)^{\alpha
-1}dy}dx}$. The quadrature approximation of $I_{\tau\times\sigma}^{i,j}$ is
denoted by $Q_{\tau\times\sigma}^{i,j}$. The associated error is
$E_{\tau\times\sigma}^{i,j}:=I_{\tau\times\sigma}^{i,j}-Q_{\tau\times\sigma
}^{i,j}$. We set
\begin{equation}
I_{\tau\times\sigma}(u,v):=\sum_{i=1}^{M}{\sum_{j=1}^{i}{u_{i}v_{j}%
I_{\tau\times\sigma}^{i,j}}}%
\end{equation}
and similarly we define $Q_{\tau\times\sigma}(u,v)$ and $E_{\tau\times\sigma
}(u,v)$. This motivates the following definition:%
\[
\mathcal{I}_{\tau}:=\{i:\tau\subset\operatorname*{supp}b_{i}\}.
\]

\begin{lemma}
Let $E_{\tau\times\sigma}^{\max}:=\max\limits_{(i,j)\in\mathcal{I}_{\tau
}\times\mathcal{I}_{\sigma},j\leqslant i}{|E_{\tau\times\sigma}^{i,j}|}$.
There exists a constant $C_{m}$ which only depends on the local polynomial
degree $m$ such that for all $u,v\in S_{\mathcal{T}}^{m}$ it holds%
\[
|E_{\tau\times\sigma}(u,v)|\leqslant C_{m}h^{-1}\Vert u\Vert_{L^{2}(\sigma
)}\Vert v\Vert_{L^{2}(\tau)}.
\]

\end{lemma}%

\proof
Let $u=\sum_{i\in\mathcal{I}_{\tau}}\alpha_{i}\left.  b_{i}\right\vert _{\tau
}$ and $v=\sum_{i\in\mathcal{I}_{\sigma}}\beta_{i}\left.  b_{i}\right\vert
_{\sigma}$. We set $%
\mbox{\boldmath$ \alpha$}%
=\left(  \alpha_{i}\right)  _{i=\mathcal{I}_{\tau}}$ and $%
\mbox{\boldmath$ \beta$}%
=\left(  \beta_{i}\right)  _{i=\mathcal{I}_{\sigma}}$ and denote by
$\left\Vert
\mbox{\boldmath$ \alpha$}%
\right\Vert $ the Euclidean norm of a vector. We estimate the error by adding
the local contributions and employing discrete Cauchy-Schwarz inequalities%
\begin{align}
|E_{\tau\times\sigma}(u,v)|  &  =\left\vert \sum_{(i,j)\in\mathcal{I}_{\tau
}\times\mathcal{I}_{\sigma},j\leqslant i}{\alpha_{i}\beta_{j}E_{\tau
\times\sigma}^{i,j}}\right\vert \\
&  \leqslant|E_{\tau\times\sigma}^{\max} |\sum_{i\in\mathcal{I}_{\tau}%
}{|\alpha_{i}|}\sum_{j\in\mathcal{I}_{\sigma}}{|\beta_{j}|}\nonumber\\
&  \leqslant\left(  m+1\right)  |E_{\tau\times\sigma}^{\max}|\Vert%
\mbox{\boldmath$ \alpha$}%
\Vert\Vert%
\mbox{\boldmath$ \beta$}%
\Vert. \label{globerror2}%
\end{align}
Well known scaling inequalities for one-dimensional finite element functions
(cf. \cite{Dahmen00}) lead to $\left\Vert
\mbox{\boldmath$ \alpha$}%
\right\Vert \leq\hat{C}_{m}h_{\tau}^{-1/2}\left\Vert u\right\Vert
_{L^{2}\left(  \tau\right)  }$ and similarly for $%
\mbox{\boldmath$ \beta$}%
$. The combination with (\ref{globerror2}) leads to the assertion with
$C_{m}=\left(  m+1\right)  C_{\operatorname*{qu}}\hat{C}_{m}^{2}$.%
\endproof

The following corollary is similar to \cite[Thm. 5.3.29]{sauter11}. We use the
fact that one-dimensional finite elements satisfy the inverse inequality%
\begin{equation}
\left\Vert u\right\Vert \leq C_{\operatorname*{inv}}h^{-s}\left\Vert
u\right\Vert _{H^{-s}\left(  \Omega\right)  }\qquad\forall u\in S_{\mathcal{T}%
}^{m}\quad\forall s\in\left[  0,\frac{1}{2}\right)  . \label{inverseinequ}%
\end{equation}
The constant $C_{\operatorname*{inv}}$ depends on $s$, $m$, $C_{\mathcal{T}}$,
and $C_{\operatorname*{qu}}$ (see \cite{Dahmen00}).

\begin{corollary}
Let $E_{K,\alpha}^{\max}:=C_{m}C_{\operatorname*{qu}}C_{\operatorname*{inv}%
}^{2}\max\limits_{\tau,\sigma\in\mathcal{T}}{E_{\tau\times\sigma}^{\max}}$ and
let $a_{K,\alpha}(\cdot,\cdot)$ be the sesquilinear form as in (\ref{auf}) and
$\tilde{a}_{K,\alpha}(\cdot,\cdot)$ the perturbed sesquilinear form in
(\ref{pervarprob}). For all $u,v\in S_{\mathcal{T}}^{m}$, it holds
\begin{equation}
|a_{K,\alpha}(u,v)-\tilde{a}_{K,\alpha}(u,v)|\leqslant E_{K,\alpha}^{\max
}h^{-2-\alpha}\Vert u\Vert_{H^{-\alpha/2}\left(  \Omega\right)  }\Vert
v\Vert_{H^{-\alpha/2}\left(  \Omega\right)  }. \label{biliconsistency}%
\end{equation}
\label{corbilin}
\end{corollary}

%

\proof
The \ quasi-uniformity of the mesh implies that $N\leq C_{\operatorname*{qu}%
}h^{-1} $. The quadrature error of the sesquilinear form can be estimated by%
\begin{align*}
&  |a_{K,\alpha}(u,v)-\tilde{a}_{K,\alpha}(u,v)|\leqslant\sum_{\tau,\sigma
\in\mathcal{T}}{|E_{\tau\times\sigma}(u,v)|}\\
&  \qquad\qquad\leqslant{C}_{m}{h^{-1}}\sum_{\tau,\sigma\in\mathcal{T}%
}{E_{\tau\times\sigma}^{\max}\Vert u\Vert_{L^{2}(\tau)}\Vert v\Vert
_{L^{2}(\sigma)}}\\
&  \qquad\qquad\leqslant C_{m}\left(  \max\limits_{\tau,\sigma\in\mathcal{T}%
}{E_{\tau\times\sigma}^{\max}}\right)  h^{-1}\sum_{\tau\in\mathcal{T}}{\Vert
u\Vert_{L^{2}(\tau)}}\sum_{\sigma\in\mathcal{T}}{\Vert v\Vert_{L^{2}(\sigma)}%
}\\
&  \qquad\qquad\leqslant C_{m}\left(  \max\limits_{\tau,\sigma\in\mathcal{T}%
}{|E_{\tau\times\sigma}^{\max}|}\right)  h^{-1}N\left\Vert u\right\Vert
\left\Vert v\right\Vert \\
&  \qquad\qquad=C_{m}C_{\operatorname*{qu}}\left(  \max\limits_{\tau,\sigma
\in\mathcal{T}}{E_{\tau\times\sigma}^{\max}}\right)  h^{-2}\Vert u\Vert\Vert
v\Vert\overset{\text{(\ref{inverseinequ})}}{\leqslant}E_{K,\alpha}^{\max
}h^{-2-\alpha}\Vert u\Vert_{H^{-\alpha/2}(\Omega)}\Vert v\Vert_{H^{-\alpha
/2}(\Omega)}.
\end{align*}%
\endproof

\begin{remark}
From Corollary \ref{corbilin} and the continuity of $a_{K,\alpha}$ it follows
that $\tilde{a}_{K,\alpha}$ is continuous in $S_{\mathcal{T}}^{m}$.
\end{remark}

%

\proof
A triangle inequality leads to%
\begin{align}
|\tilde{a}_{K,\alpha}(u,v)|  &  \leqslant|a_{K,\alpha}(u,v)|+|\tilde
{a}_{K,\alpha}(u,v)-a_{K,\alpha}(u,v)|\label{defCchat}\\
&  \leqslant\hat{C}_{c}\Vert u\Vert_{H^{-\alpha/2}\left(  \Omega\right)
}\Vert v\Vert_{H^{-\alpha/2}\left(  \Omega\right)  }\qquad\text{with }\hat
{C}_{c}:=E_{K,\alpha}^{\max}h^{-2-\alpha}+C_{c}.\nonumber
\end{align}%
\endproof

For the functional $G(v)$ as in (\ref{discprob}), $v\in S_{\mathcal{T}}^{m}$,
$m\in\mathbb{N}_{0}$, it holds%
\begin{equation}
G(v)=\int_{0}^{1}{gv}=\sum_{i=1}^{M}{\sum_{\tau\subset\operatorname*{supp}%
(b_{i})}{\beta_{i}}}\int_{\tau}gb_{i}. \label{luf}%
\end{equation}
We proceed similar as for $a(\cdot,\cdot)$ and define $I_{\tau}^{i}%
:=\int_{\tau}{gb_{i}}$ and the quadrature operator $Q_{\tau}^{i}$ in the same
fashion so that the error is $E_{\tau}^{i}=I_{\tau}^{i}-Q_{\tau}^{i}$. We set
$I_{\tau}(v):=\sum_{i\in\mathcal{I}_{\tau}}{\beta_{i}I_{\tau}^{i}}$ and define
$Q_{\tau}(v)$ and $E_{\tau}(v)$ analogously. With $E_{\tau}^{\max}=\max
_{i\in\mathcal{I}_{\tau}}{|E_{\tau}^{i}|}$ we get%
\[
|E_{\tau}(v)|=|\sum_{i\in\mathcal{I}_{\tau}}{\beta_{i}E_{\tau}^{i}}%
|\leqslant\sqrt{m+1}E_{\tau}^{\max}\Vert%
\mbox{\boldmath$ \beta$}%
\Vert\leqslant E_{\tau}^{\max}\hat{C}_{m}\sqrt{\frac{m+1}{h_{\tau}}}\Vert
v\Vert_{L^{2}\left(  \tau\right)  }.
\]

\begin{corollary}
\label{CorGdiff}Let $E_{G}^{\max}:={\sqrt{C_{m}C_{\operatorname*{qu}}}%
C}_{\operatorname*{inv}}\max\limits_{\tau\in\mathcal{T}}{E_{\tau}^{\max}}$.
Let $G(\cdot)$ be the functional in (\ref{luf}) and $\tilde{G}(\cdot)$ the
perturbed functional as in (\ref{pervarprob}). For all $v\in S_{\mathcal{T}%
}^{m}$ it holds%
\begin{equation}
|G(v)-\tilde{G}(v)|\leqslant E_{G}^{\max}h^{-(1+\alpha)/2}\Vert v\Vert
_{H^{-\alpha/2}\left(  \Omega\right)  }. \label{difffuncG}%
\end{equation}
\label{functional}
\end{corollary}

%

\proof
Again we localize the difference $G-\tilde{G}$ and obtain%
\begin{align*}
|G(v)-\tilde{G}(v)|  &  \leqslant\sum_{\tau\in\mathcal{T}}{|E_{\tau}%
(v)|}\leqslant\sum_{\tau\in\mathcal{T}}{E_{\tau}^{\max}\sqrt{C_{m}}%
h^{-1/2}\Vert v\Vert_{L^{2}(\tau)}}\\
&  \leqslant{\sqrt{C_{m}}h^{-1/2}}\left(  \max_{\tau\in\mathcal{T}}{E_{\tau
}^{\max}}\right)  \sum_{\tau\in\mathcal{T}}{\Vert v\Vert_{L^{2}(\tau)}}\\
&  \leqslant\sqrt{C_{m}Q_{\operatorname*{qu}}}\left(  \max_{\tau\in
\mathcal{T}}{E_{\tau}^{\max}}\right)  h^{-1}\Vert v\Vert.
\end{align*}
An inverse inequality for $v$ leads to the assertion.%
\endproof

We will now investigate the error of the perturbed Galerkin method. The proof
of the next theorem is based on the first Strang lemma.

\begin{theorem}
Suppose that the kernel function $K$ in Abel's integral equation is
$s$-admissible for $s=-\alpha/2$. Assume that the quadrature method is
sufficiently accurate such that%
\begin{equation}
E_{K,\alpha}^{\max}h^{-2-\alpha}\leq\frac{\tilde{\gamma}}{2} \label{condquad}%
\end{equation}
with $\tilde{\gamma}$ as in (\ref{defgammatildeconst}). Then, the perturbed
Galerkin method (\ref{pervarprob}) is $H^{-\alpha/2}\left(  \Omega\right)
$-elliptic and the fully discrete equations (\ref{pervarprob}) have a unique
solution $\tilde{f}_{S}$ that satisfies the error estimate
\begin{equation}
\Vert f-\tilde{f}_{S}\Vert_{H^{-\alpha/2}\left(  \Omega\right)  }\leqslant
C\left(  \min_{w\in S}\left(  \Vert f-w\Vert_{H^{-\alpha/2}\left(
\Omega\right)  }+E_{K,\alpha}^{\max}h^{-2-\alpha}\Vert w\Vert_{H^{-\alpha
/2}\left(  \Omega\right)  }\right)  +E_{G}^{\max}h^{-(1+\alpha)/2}\right)  ,
\label{errorperturbed}%
\end{equation}
for some $C>0$.
\end{theorem}

%

\proof
Theorem \ref{main1} implies that the exact sesquilinear form $a_{K,\alpha}$ is
$H^{-\alpha/2}\left(  \Omega\right)  $-elliptic and we derive the property for
the perturbed version $\tilde{a}_{K,\alpha}$ next. From Corollary
\ref{corbilin} we conclude that
\[
|a_{K,\alpha}(u,v)-\tilde{a}_{K,\alpha}(u,v)|\leqslant E_{K,\alpha}^{\max
}h^{-2-\alpha}\Vert u\Vert_{H^{-\alpha/2}\left(  \Omega\right)  }\Vert
v\Vert_{H^{-\alpha/2}\left(  \Omega\right)  }.
\]
This implies%
\begin{align*}
\operatorname{Re}\tilde{a}_{K,\alpha}(u,u)  &  \geqslant\operatorname{Re}%
a_{K,\alpha}(u,u)-\operatorname{Re}(a_{K,\alpha}(u,u)-\tilde{a}_{K,\alpha
}(u,u))\\
&  \geqslant\tilde{\gamma}\Vert u\Vert_{H^{-\alpha/2}\left(  \Omega\right)
}^{2}-|a_{K,\alpha}(u,u)-\tilde{a}_{K,\alpha}(u,u)|\\
&  \geqslant(\tilde{\gamma}-E_{K,\alpha}^{\max}h^{-2-\alpha})\Vert
u\Vert_{H^{-\alpha/2}\left(  \Omega\right)  }^{2}.
\end{align*}
Hence, condition (\ref{condquad}) implies the $H^{-\alpha/2}\left(
\Omega\right)  $-ellipticity of the perturbed sesquilinear form. From the
Lax-Milgram lemma we conclude that (\ref{pervarprob}) has a unique solution
$\tilde{f}_{S}$. Next we prove the error estimate (\ref{errorperturbed})
\begin{align}
\Vert f-\tilde{f}_{S}\Vert_{H^{-\alpha/2}\left(  \Omega\right)  }  &
\leqslant\Vert f-f_{S}\Vert_{H^{-\alpha/2}\left(  \Omega\right)  }+\Vert
f_{S}-\tilde{f}_{S}\Vert_{H^{-\alpha/2}\left(  \Omega\right)  } \label{quader}%
\\
&  \leqslant\Vert f-f_{S}\Vert_{H^{-\alpha/2}\left(  \Omega\right)  }+\frac
{2}{\tilde{\gamma}}\frac{\operatorname{Re}\tilde{a}_{K,\alpha}(f_{S}-\tilde
{f}_{S},f_{S}-\tilde{f}_{S})}{\Vert f_{S}-\tilde{f}_{S}\Vert_{H^{-\alpha
/2}\left(  \Omega\right)  }}\nonumber\\
&  \leqslant\Vert f-f_{S}\Vert_{H^{-\alpha/2}\left(  \Omega\right)  }+\frac
{2}{\tilde{\gamma}}\sup_{v\in S\backslash\{0\}}\frac{|\tilde{a}_{K,\alpha
}(f_{S}-\tilde{f}_{S},v)|}{\Vert v\Vert_{H^{-\alpha/2}\left(  \Omega\right)
}}\nonumber\\
&  =\Vert f-f_{S}\Vert_{H^{-\alpha/2}\left(  \Omega\right)  }+\frac{2}%
{\tilde{\gamma}}\sup_{v\in S\backslash\{0\}}\frac{|\tilde{a}_{K,\alpha}%
(f_{S},v)-\tilde{G}(v)|}{\Vert v\Vert_{H^{-\alpha/2}\left(  \Omega\right)  }%
}\nonumber\\
&  \leqslant\Vert f-f_{S}\Vert_{H^{-\alpha/2}\left(  \Omega\right)  }+\frac
{2}{\tilde{\gamma}}\sup_{v\in S\backslash\{0\}}\frac{|\tilde{a}_{K,\alpha
}(f_{S},v)-a_{K,\alpha}(f_{S},v)|+|G(v)-\tilde{G}(v)|}{\Vert v\Vert
_{H^{-\alpha/2}\left(  \Omega\right)  }}.\nonumber
\end{align}
We consider the difference $|\tilde{a}_{K,\alpha}(f_{S},v)-a_{K,\alpha}%
(f_{S},v)|$ and obtain by using the continuity of $a_{K,\alpha}$ and
$\tilde{a}_{K,\alpha}$ as well as the consistency (\ref{biliconsistency})
that
\begin{align*}
&  |\tilde{a}_{K,\alpha}(f_{S},v)-a_{K,\alpha}(f_{S},v)|\leqslant|\tilde
{a}_{K,\alpha}(f_{S}-w,v)|+|\tilde{a}_{K,\alpha}(w,v)-a_{K,\alpha
}(w,v)|+|a_{K,\alpha}(w-f_{S},v)|\\
&  \qquad\qquad\leqslant\left(  \tilde{C}_{c}+\frac{\tilde{\gamma}}{2}\right)
\Vert f_{S}-w\Vert_{H^{-\alpha/2}\left(  \Omega\right)  }\Vert v\Vert
_{H^{-\alpha/2}\left(  \Omega\right)  }+E_{K,\alpha}^{\max}h^{-2-\alpha}\Vert
w\Vert_{H^{-\alpha/2}\left(  \Omega\right)  }\Vert v\Vert_{H^{-\alpha
/2}\left(  \Omega\right)  }\\
&  \qquad\qquad\quad+\tilde{C}_{c}\Vert w-f_{S}\Vert_{H^{-\alpha/2}\left(
\Omega\right)  }\Vert v\Vert_{H^{-\alpha/2}\left(  \Omega\right)  }%
\end{align*}
holds for all $w\in S$. This leads to
\[
\sup_{v\in S\backslash\{0\}}\frac{|\tilde{a}(f_{S},v)-a(f_{S},v)|}{\Vert
v\Vert_{H^{-\alpha/2}\left(  \Omega\right)  }}\leqslant\min_{w\in S}\left(
C\Vert f_{S}-w\Vert_{H^{-\alpha/2}\left(  \Omega\right)  }+E_{K,\alpha}^{\max
}h^{-2-\alpha}\Vert w\Vert_{H^{-\alpha/2}\left(  \Omega\right)  }\right)
\]
with $C:=2\tilde{C}_{c}+\frac{\tilde{\gamma}}{2}$. The combination with
(\ref{discerror1}) and (\ref{difffuncG}) leads to%
\begin{align*}
\Vert f-\tilde{f}_{S}\Vert_{H^{-\alpha/2}(\Omega)} \leqslant &  C \min_{w\in
S}\Big( \Vert f-w\Vert_{H^{-\alpha/2}(\Omega)}\\
&  +\frac{2}{\tilde{\gamma}}(\Vert f-w\Vert_{H^{-\alpha/2}(\Omega)
}+E_{K,\alpha}^{\max}h^{-2-\alpha}\Vert w\Vert_{H^{-\alpha/2}(\Omega)}%
+\sup_{v\in S\backslash\{0\}}\frac{|G(v)-\tilde{G}(v)|}{\Vert v\Vert
_{H^{-\alpha/2}(\Omega)}})\Big)\\
\leqslant &  C \Big( \min_{w\in S}\big(\Vert f-w \Vert_{H^{-\alpha/2}(\Omega
)}+E_{K,\alpha}^{\max}h^{-2-\alpha}\Vert w\Vert_{H^{-\alpha/2}(\Omega
)}\big)+E_{G}^{\max}h^{-(1+\alpha)/2}\Big).
\end{align*}
\endproof

\subsection{Choice of the Quadrature Order\label{SecQuadOrder}}

The results of Sections \ref{SubSecPert} and \ref{SecInfluence} allow to
determine the appropriate quadrature order $n$ so that the perturbed Galerkin
method converges at the same rate as the unperturbed Galerkin method.

\begin{theorem}
Let the assumption of Theorem \ref{main1} be satisfied. Assume that the kernel
function $K$ satisfies (\ref{analyticK}) and that the right-hand side $g$
fulfills condition (\ref{ganalytic}). Let the quadrature order in
(\ref{gauserrgal}) for the integrals over those $\tau\times\sigma$ which
satisfy
\[
\operatorname*{dist}\left(  \tau,\sigma\right)  \geq\Lambda_{K}\max\left\{
h_{\tau},h_{\sigma}\right\}
\]
be chosen according to
\begin{equation}
n_{1}=\left\lceil \left(  m+2+\alpha/4\right)  \frac{\log\left(  \frac{1}%
{h}\right)  }{\log\left(  \frac{2}{\Lambda_{K}}\frac{\operatorname*{dist}%
\left(  \tau,\sigma\right)  }{\max\left\{  h_{\tau},h_{\sigma}\right\}
}\right)  }\right\rceil . \label{orderni}%
\end{equation}
Assume that the remaining singular and near-singular integrals are evaluated
exactly. Then the fully discrete Galerkin method is stable. Let the number of
quadrature points for the right-hand side be chosen as
\[
n_{2}=\left\lceil m+\frac{\alpha}{2}+\frac{3}{4}\right\rceil .
\]
If the exact solution is $H^{m+1}(\Omega)$, then, the perturbed Galerkin
method converges at the same rate as the Galerkin method.
\label{errorestimate}
\end{theorem}

\begin{remark}
The assumption that all singular integrals are evaluated exactly is related to
the fact that we have omitted the error analysis for the singular integrals
(cf. Remark \ref{RemExact}). By using techniques as developed in
\cite{sauter11} this case can be handled while the technicalities are
increased. We emphasize that the evaluation of the \emph{nearly singular} and
\emph{regular} integrals is dominating the cost for building the system matrix
since their number is $O\left(  M^{2}\right)  $ while the number of
\emph{singular} integrals is $O\left(  M\right)  $. Quadrature methods for the
singular integrals will be presented in a forthcoming paper.
\end{remark}

\proof
We employ Lemma \ref{TheoLocQuad} and standard inverse estimates for finite
element functions, i.e., there exists a constant $C_{\operatorname*{inv}}$
such that $\Vert b_{i}\Vert_{C^{m}(\tau)}\leq C_{\operatorname*{inv}}h_{\tau
}^{-m}$ (we may use the same notation as in (\ref{inverseinequ}) by a suitable
adjustment of $C_{\operatorname*{inv}}$). Thus, we obtain%
\[
E_{\tau\times\sigma}^{\max}\leq2C_{m}C_{K}C_{\operatorname*{inv}%
}\operatorname*{e}\Lambda_{K}h_{\tau}^{-m}\max_{\tau,\sigma\in\mathcal{T}%
}\left(  \frac{1}{\operatorname*{dist}\nolimits^{1-\alpha}\left(  \tau
,\sigma\right)  }\left(  \frac{\Lambda_{K}}{2}\frac{\max\left\{  h_{\tau
},h_{\sigma}\right\}  }{\operatorname*{dist}\left(  \tau,\sigma\right)
}\right)  ^{2n_{1}}\right)  .
\]
The assumption $\operatorname*{dist}\left(  \tau,\sigma\right)  \geq
\Lambda_{K}\max\left\{  h_{\tau},h_{\sigma}\right\}  \geq\Lambda
_{K}C_{\operatorname*{qu}}^{-1}h$ implies%
\[
E_{\tau\times\sigma}^{\max}\leq2C_{m}C_{K}C_{\operatorname*{inv}%
}C_{\operatorname*{qu}}^{1-\alpha}\operatorname*{e}\Lambda_{K}^{\alpha
}h^{\alpha-1-m}\left(  \frac{\Lambda_{K}}{2}\frac{\max\left\{  h_{\tau
},h_{\sigma}\right\}  }{\operatorname*{dist}\left(  \tau,\sigma\right)
}\right)  ^{2n_{1}}.
\]
From (\ref{condquad}) we deduce that the fully discrete Galerkin method is
stable if
\[
2C_{m}C_{K}C_{\operatorname*{inv}}C_{\operatorname*{qu}}^{1-\alpha
}\operatorname*{e}\Lambda_{K}^{\alpha}h^{-3-m}\left(  \frac{\Lambda_{K}}%
{2}\frac{\max\left\{  h_{\tau},h_{\sigma}\right\}  }{\operatorname*{dist}%
\left(  \tau,\sigma\right)  }\right)  ^{2n_{1}}\leq\frac{\tilde{\gamma}}{2}.
\]
To obtain an optimal convergence order we obtain from (\ref{errorperturbed})
the stronger condition%
\[
2C_{m}C_{K}C_{\operatorname*{inv}}C_{\operatorname*{qu}}^{1-\alpha
}\operatorname*{e}\Lambda_{K}^{\alpha}h^{-3-m}\left(  \frac{\Lambda_{K}}%
{2}\frac{\max\left\{  h_{\tau},h_{\sigma}\right\}  }{\operatorname*{dist}%
\left(  \tau,\sigma\right)  }\right)  ^{2n_{1}}\leq Ch^{m+1+\alpha/2}.
\]
We solve this last condition for $n_{1}$ and get%
\begin{equation}
n_{1}=\left\lceil \left(  m+2+\alpha/4\right)  \frac{\log\left(  \frac{1}%
{h}\right)  }{\log\left(  \frac{2}{\Lambda_{K}}\frac{\operatorname*{dist}%
\left(  \tau,\sigma\right)  }{\max\left\{  h_{\tau},h_{\sigma}\right\}
}\right)  }\right\rceil . \label{n1quadorder}%
\end{equation}
Finally, we have to determine the number $n_{2}$ of Gauss points for the
approximation of the right-hand side. To preserve the optimal convergence
order, the last term in (\ref{errorperturbed}) has to be bounded by
$Ch^{m+1+\alpha/2}$. The combination with Corollary \ref{CorGdiff} leads to
the condition
\[
{\sqrt{C_{m}C_{\operatorname*{qu}}}C}_{\operatorname*{inv}}\left(
\max\limits_{\tau\in\mathcal{T}}{E_{\tau}^{\max}}\right)  h^{-(1+\alpha
)/2}\leq Ch^{m+1+\alpha/2}.
\]
The quantity ${E_{\tau}^{\max}}$ can be estimated via (\ref{gausserrrhs}) so
that%
\[
{\sqrt{C_{m}C_{\operatorname*{qu}}}C}_{\operatorname*{inv}}^{2}2\Lambda
_{g}\left(  \frac{\Lambda_{g}}{2}h_{\tau}\right)  ^{2n_{2}}h^{-(1+\alpha
)/2-m}\leq Ch^{m+1+\alpha/2}%
\]
is a sufficient condition. We solve this for $n_{2}$ and obtain%
\[
n_{2}\leq\left(  \frac{2m+3/2+\alpha}{2}\right)  \frac{\log\frac{1}{h}}%
{\log\frac{2}{\Lambda_{g}}+\log\frac{1}{h_{\tau}}}\leq\frac{2m+3/2+\alpha}%
{2}.
\]
\endproof

\section{Numerical Experiments\label{SecNumexp}}

The numerical experiments are used to investigate the sharpness of our
theoretical error estimates. For this we study how the quadrature method
influences the convergence rate. The first test is done for the integral
equation
\begin{equation}
\Gamma({1}/{2})^{-1}\int_{0}^{x}\frac{K\left(  x,y\right)  }{\left(
x-y\right)  ^{1/2}}f(y)\mathrm{d}y=g(x)\quad\forall x\in\Omega\label{NumExp1}%
\end{equation}
for $g\left(  x\right)  :=-\frac{\sqrt{\pi}x^{2}}{160}\left(  -60+x\left(
11+5x\right)  \right)  $ and $K(x,y)=1-\frac{x+y}{10}-\frac{xy}{10}$. Note
that the kernel function has the representation as in (\ref{repK1}) for
$\psi_{2}\left(  x\right)  =1+x$ and $\psi_{m}=1$ otherwise. The coefficients
$d_{n,m}$ then equal $d_{1,1}=\frac{11}{10}$, $d_{2,2}=-\frac{1}{10}$, and
$d_{n,m}=0$ otherwise. Clearly the constants $c_{n,s}$, $C_{n,s}$ for
$n\in\mathbb{N}\backslash\left\{  2\right\}  $ equal $1$. For $n=2$ we obtain
from \cite[Lemma 2.0.3]{HeimgartnerMaster} for $s=-\alpha/2\in\left(
-1/2,0\right)  $ that the choices $C_{2,s}=2\sqrt{2}$ and $c_{2,s}=2^{-1/2}$
are admissible in (\ref{cnsCns}) (for $n=2$) and
\begin{align*}
\sum_{n,m=1}^{\infty}|d_{n,m}|C_{n,s}C_{m,s}  &  =\frac{11}{10}C_{1,s}%
^{2}+\frac{1}{10}C_{2,s}^{2}=\frac{19}{10}=:C_{s}^{2}\\
\gamma\sum_{n=1}^{\infty}d_{n,n}c_{n,s}^{2}-C_{c}\sum_{\substack{n,m=1\\n\neq
m}}^{\infty}|d_{n,m}|C_{n,s}C_{m,s}  &  =\gamma\left(  \frac{11}{10}%
c_{1,s}^{2}-\frac{1}{10}c_{2,s}^{2}\right)  =\frac{21}{20}\gamma>0.
\end{align*}
Thus, $K\in\mathcal{A}\left(  -\frac{\alpha}{2}\right)  $. The solution is
given by $f(y)=y^{3/2}$. For the discretization in this example, we have used
a uniform mesh with width $h=N^{-1}$.

\subsection{Changing the Quadrature Order}

Let $s_{i}:=m+i+\alpha/4$ so that the prefactor in (\ref{orderni}) in front of
the logarithmic terms equals $s_{2}$. We have varied $i\in\left\{
1,2,\ldots,5\right\}  $ to numerically validate the sharpness of estimate
(\ref{orderni}) for the number of quadrature points.

\begin{figure}[ptb]
\begin{tikzpicture}
[scale=0.8, transform shape]
\begin{axis}[title=Error of Galerkin Method on $S_\mathcal{T}^1$,
xlabel={Partition Size $N$},xmode=log, log basis x={2},ylabel={$\Vert f-\tilde{f}_S\Vert_{H^{-1/4}}$}, ymode=log, legend entries={$s_1$,$s_2$,$s_3$,$s_4$,$s_5$}, legend pos=south west]
\addplot coordinates {
(2^5,2.11E-04)
(2^6,6.34E-05)
(2^7, 1.93E-05)
(2^8,5.75E-06)
(2^9,1.71E-06)
(2^10,5.02E-07)
(2^11,1.43E-07)
(2^12,4.00E-08)
};
\addplot coordinates {
(2^5,1.79E-04)
(2^6,5.05E-05)
(2^7, 1.43E-05)
(2^8,4.08E-06)
(2^9,1.17E-06)
(2^10,3.39E-07)
(2^11,9.88E-08)
(2^12,2.89E-08)
};
\addplot coordinates {
(2^5,1.79E-04)
(2^6,5.05E-05)
(2^7,1.43E-05 )
(2^8,4.08E-06)
(2^9,1.17E-06)
(2^10,3.39E-07)
(2^11,9.88E-08)
(2^12,2.89E-08)
};
\addplot coordinates {
(2^5,1.79E-04)
(2^6,5.05E-05)
(2^7,1.43E-05)
(2^8,4.08E-06)
(2^9,1.17E-06)
(2^10,3.39E-07)
(2^11,9.88E-08)
(2^12,2.89E-08)
};
\addplot coordinates {
(2^5,1.79E-04)
(2^6,5.05E-05)
(2^7, 1.43E-05)
(2^8,4.08E-06)
(2^9,1.17E-06)
(2^10,3.39E-07)
(2^11,9.88E-08)
(2^12,2.89E-08)
};
\end{axis}
\draw[red,line width=0.7mm, dashed](4,4) -- (5.1,3);
\end{tikzpicture}
\begin{tikzpicture}
[scale=0.8, transform shape]
\begin{axis}[title=Error of Galerkin Method on $S_\mathcal{T}^0$, xlabel={Partition Size $N$},xmode=log, log basis x={2},ylabel={$\Vert f-\tilde{f}_S\Vert_{H^{-1/4}}$}, ymode=log, legend entries={$s_1$,$s_2$,$s_3$,$s_4$,$s_5$}, legend pos=south west]
\addplot coordinates {
(2^5,0.01278108)
(2^6,0.006350216)
(2^7,0.003161266)
(2^8,0.001580652)
(2^9,0.000790265)
(2^10,3.96E-04)
(2^11,1.98E-04)
(2^12,9.91E-05)
};
\addplot coordinates {
(2^5,0.012910279)
(2^6,0.006412444)
(2^7,0.003195482)
(2^8,0.001595044)
(2^9,0.000796845)
(2^10,3.98E-04)
(2^11,1.99E-04)
(2^12,9.95E-05)
};
\addplot coordinates {
(2^5,0.012911173)
(2^6,0.006412648)
(2^7,0.003195527)
(2^8,0.001595053)
(2^9,0.000796847)
(2^10,3.98E-04)
(2^11,1.99E-04)
(2^12,9.95E-05)
};
\addplot coordinates {
(2^5,0.012911185)
(2^6,0.006412651)
(2^7,0.003195527)
(2^8,0.001595053)
(2^9,0.000796847)
(2^10,3.98E-04)
(2^11,1.99E-04)
(2^12,9.95E-05)
};
\addplot coordinates {
(2^5,0.012911185)
(2^6,0.006412651)
(2^7,0.003195527)
(2^8,0.001595053)
(2^9,0.000796847)
(2^10,3.98E-04)
(2^11,1.99E-04)
(2^12,9.95E-05)
};
\end{axis}
\draw[red,line width = 0.7mm, dashed](4,4) -- (5.1,3);;
\end{tikzpicture}
\caption{Plots of the error in logarithmic scale and the dashed reference line
illustrates, on the left, quadratic convergence and, on the right, linear
convergence. }%
\label{figs1}%
\end{figure}
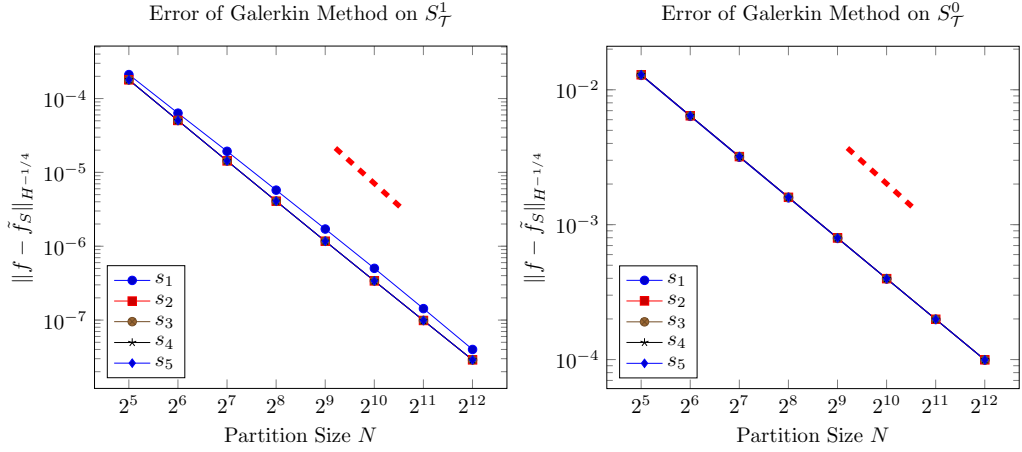

Figure \ref{figs1} shows that, for properly chosen numbers of quadrature
points (cf. Sec. \ref{SecQuadOrder}), the slopes of the convergence curves
match with the theoretical predicted linear and quadratic convergence (dashed
lines).

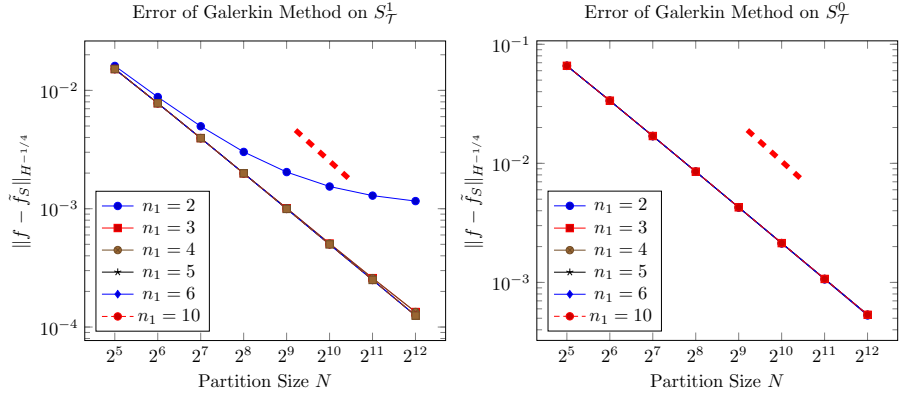
\begin{figure}[ptb]
\begin{tikzpicture}[scale=0.7, transform shape]
\begin{axis}[title=Error of Galerkin Method on $S_\mathcal{T}^1$, xlabel={Partition Size $N$},xmode=log, log basis x={2},ylabel={$\Vert f-\tilde{f}_S\Vert_{H^{-1/4}}$}, ymode=log, legend entries={$n_{1}=2$,$n_{1}=3$,$n_{1}=4$,$n_{1}=5$,$n_{1}=6$,$n_{1}=10$}, legend pos=south west]
\addplot coordinates {
(2^5,0.016090655)
(2^6,0.008786697)
(2^7,0.004972899)
(2^8,3.02E-03)
(2^9,2.04E-03)
(2^10, 1.54E-03)
(2^11,1.29E-03)
(2^12,1.16E-03)
};
\addplot coordinates {
(2^5,0.015151019)
(2^6,0.00779146)
(2^7,0.003953966)
(2^8,1.99E-03)
(2^9,1.01E-03)
(2^10,5.08E-04)
(2^11,2.59E-04)
(2^12,1.34E-04)
};
\addplot coordinates {
(2^5,0.015151019)
(2^6,0.00779146)
(2^7,0.003953966)
(2^8,1.99E-03)
(2^9,1.01E-03)
(2^10,5.08E-04)
(2^11,2.59E-04)
(2^12,1.34E-04)
};
\addplot coordinates {
(2^5,0.015142855)
(2^6,0.007782965)
(2^7,0.003945345)
(2^8,1.99E-03)
(2^9,9.97E-04)
(2^10,4.99E-04)
(2^11,2.50E-04)
(2^12,1.25E-04)
};
\addplot coordinates {
(2^5,0.01514275)
(2^6,0.007782856)
(2^7,0.003945234)
(2^8,1.99E-03)
(2^9,9.97E-04)
(2^10,4.99E-04)
(2^11,2.50E-04)
(2^12,1.25E-04)
};
\addplot coordinates {
(2^5,0.015142748)
(2^6,0.007782854)
(2^7,0.003945232)
(2^8,1.99E-03)
(2^9,9.96E-04)
(2^10,4.99E-04)
(2^11,2.50E-04)
(2^12,1.25E-04)
};
\addplot coordinates {
(2^5,0.015142748)
(2^6,0.007782854)
(2^7,0.003945232)
(2^8,1.99E-03)
(2^9,9.96E-04)
(2^10,4.99E-04)
(2^11,2.50E-04)
(2^12,1.25E-04)
};
\end{axis}
\draw[red,line width=0.7mm, dashed](4,4) -- (5.1,3);
\end{tikzpicture}
\begin{tikzpicture}[scale=0.7, transform shape]
\begin{axis}[title=Error of Galerkin Method on $S_\mathcal{T}^0$, xlabel={Partition Size $N$},xmode=log, log basis x={2},ylabel={$\Vert f-\tilde{f}_S\Vert_{H^{-1/4}}$}, ymode=log, legend entries={$n_{1}=2$,$n_{1}=3$,$n_{1}=4$,$n_{1}=5$,$n_{1}=6$,$n_{1}=10$}, legend pos=south west]
\addplot coordinates {
(2^5,0.065885828)
(2^6,0.033553456)
(2^7,0.016925517)
(2^8,0.008495077)
(2^9,0.004251892)
(2^10,0.002124367)
(2^11,0.0010599)
(2^12,0.000528043)
};
\addplot coordinates {
(2^5,0.065956371)
(2^6,0.033606141)
(2^7,0.01696377)
(2^8,0.00852248)
(2^9,0.004271395)
(2^10,0.002138202)
(2^11,0.001069699)
(2^12,0.000534977)
};
\addplot coordinates {
(2^5,0.065957497)
(2^6,0.033606979)
(2^7,0.016964378)
(2^8,0.008522916)
(2^9,0.004271705)
(2^10,0.002138422)
(2^11,0.001069854)
(2^12,0.000535088)
};
\addplot coordinates {
(2^5,0.065957518)
(2^6,0.033606995)
(2^7,0.01696439)
(2^8,0.008522924)
(2^9,0.004271711)
(2^10,0.002138426)
(2^11,0.001069857)
(2^12,0.00053509)
};
\addplot coordinates {
(2^5,0.065957518)
(2^6,0.033606995)
(2^7,0.01696439)
(2^8,0.008522924)
(2^9,0.004271711)
(2^10,0.002138427)
(2^11,0.001069857)
(2^12,0.00053509)
};
\addplot coordinates {
(2^5,0.065957518)
(2^6,0.033606995)
(2^7,0.01696439)
(2^8,0.008522924)
(2^9,0.004271711)
(2^10,0.002138427)
(2^11,0.001069857)
(2^12,0.00053509)
};
\end{axis}
\draw[red,line width = 0.7mm, dashed](4,4) -- (5.1,3);
\end{tikzpicture}
\caption{Plots of the error in logarithmic scale and the dashed reference line
illustrates, on the left, quadratic convergence and, on the right, linear
convergence. }%
\label{figm}%
\end{figure}

\textbf{Fixed Quadrature Order:} In the previous experiment we adapted the
quadrature order depending on $\operatorname*{dist}(\tau,\sigma)$. Now, we fix
a quadrature order for the whole experiment. The results are depicted in
Figure \ref{figm} and clearly show that a low quadrature order for
$S_{\mathcal{T}}^{1}$ significantly pollutes the convergence rate.

\subsection{Dependence of the Solution of Abel's Integral Equation on the
Order of Singularity $\alpha$}

In this experiment, we investigate the sensitivity of our method on the
parameter $\alpha\in\left(  0,1\right)  $ in Abel's integral equation. We
consider the equation%

\begin{equation}
\frac{1}{\Gamma(\alpha)}\int_{0}^{x}(x-y)^{\alpha-1}K(x,y)f(y)\mathrm{d}%
y=g(x)\quad\forall x\in\Omega\label{numexp1}%
\end{equation}
for $K\left(  x,y\right)  $ as in (\ref{NumExp1}) and $g\left(  x\right)
:=\frac{\left(  2-\alpha\right)  \left(  1-\alpha\right)  \pi x^{2}}%
{\Gamma\left(  \alpha\right)  \sin\left(  \alpha\pi\right)  }\left(
30-x\left(  6-\alpha+x\left(  3-\alpha\right)  \right)  \right)  $. The
solution is given by $f_{\alpha}(y)=y^{2-\alpha}$. We choose $\alpha=\frac
{k}{10}$, for $k\in\{1,\ldots,9\}$, and investigate the error in Figure
\ref{alphachange}. 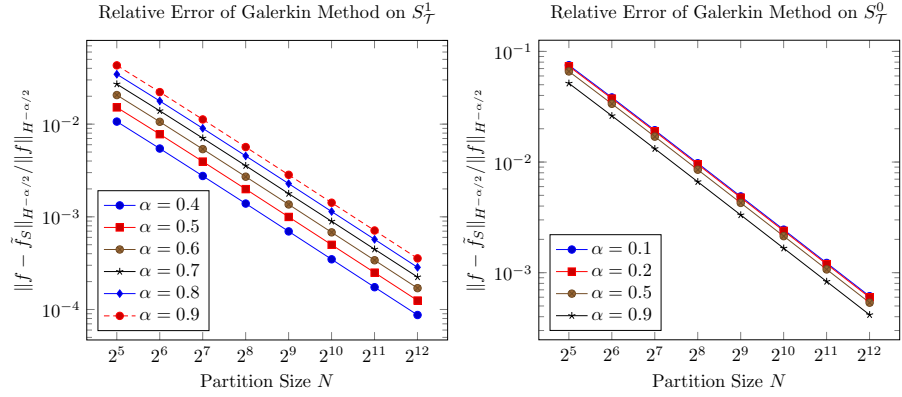
\begin{figure}[ptb]
\begin{tikzpicture}[scale=0.7, transform shape]
\begin{axis}[title=Relative Error of Galerkin Method on $S_\mathcal{T}^1$, xlabel={Partition Size $N$},xmode=log, log basis x={2},ylabel={$\Vert f-\tilde{f}_S\Vert_{H^{-\alpha/2}}/\Vert f \Vert_{H^{-\alpha/2}}$}, ymode=log, legend entries={ $\alpha=0.4$,$\alpha=0.5$,$\alpha=0.6$,$\alpha=0.7$,$\alpha=0.8$,$\alpha=0.9$}, legend pos=south west]
\addplot coordinates {
(2^5, 0.01065488)
(2^6, 0.005454185)
(2^7,0.002758987)
(2^8,1.39E-03)
(2^9,6.96E-04)
(2^10,3.48E-04)
(2^11,1.74E-04)
(2^12,8.72E-05)
};
\addplot coordinates {
(2^5,0.015154843 )
(2^6, 0.007788319)
(2^7,0.003946882)
(2^8,1.99E-03)
(2^9,9.97E-04)
(2^10,4.99E-04)
(2^11,2.50E-04)
(2^12,1.25E-04)
};
\addplot coordinates {
(2^5, 0.020593439)
(2^6, 0.010598825)
(2^7,0.005375841)
(2^8,2.71E-03)
(2^9,1.36E-03)
(2^10,6.80E-04)
(2^11,3.40E-04)
(2^12,1.70E-04)
};
\addplot coordinates {
(2^5, 0.02703713)
(2^6, 0.013919483)
(2^7,0.007060704)
(2^8,3.56E-03)
(2^9,1.78E-03)
(2^10,8.94E-04)
(2^11,4.47E-04)
(2^12,2.24E-04)
};
\addplot coordinates {
(2^5,0.034530474 )
(2^6,0.017773283 )
(2^7,0.009013297)
(2^8,4.54E-03)
(2^9,2.28E-03)
(2^10,1.14E-03)
(2^11,5.71E-04)
(2^12,2.85E-04)
};
\addplot coordinates {
(2^5,0.043065432)
(2^6, 0.022172663)
(2^7,0.011240757)
(2^8,5.66E-03)
(2^9,2.84E-03)
(2^10,1.42E-03)
(2^11,7.12E-04)
(2^12,3.56E-04)
};
\end{axis}
\end{tikzpicture}
\begin{tikzpicture}[scale=0.7, transform shape]
\begin{axis}[title=Relative Error of Galerkin Method on $S_\mathcal{T}^0$, xlabel={Partition Size $N$},xmode=log, log basis x={2},ylabel={$\Vert f-\tilde{f}_S\Vert_{H^{-\alpha/2}}/\Vert f \Vert_{H^{-\alpha/2}}$}, ymode=log, legend entries={$\alpha=0.1$,$\alpha=0.2$, $\alpha=0.5$,$\alpha=0.9$}, legend pos=south west]
\addplot coordinates {
(2^5,0.075089887)
(2^6,0.038371763)
(2^7,0.019397371)
(2^8,0.009752244)
(2^9,0.004889578)
(2^10,0.002448164)
(2^11,0.001224926)
(2^12,0.000612674)
};
\addplot coordinates {
(2^5,0.073350941)
(2^6,0.037448248)
(2^7,0.018922037)
(2^8,0.009511113)
(2^9,0.004768172)
(2^10,0.002387249)
(2^11,0.001194417)
(2^12,0.000597407
)
};
\addplot coordinates {
(2^5,0.065937769)
(2^6,0.033602461)
(2^7,0.01696338)
(2^8,0.008522723)
(2^9,0.004271671)
(2^10, 0.002138419)
(2^11,0.001069856)
(2^12,0.00053509)
};
\addplot coordinates {
(2^5,0.051464024)
(2^6,0.026184296)
(2^7,0.01320181)
(2^8,0.006628599)
(2^9,0.003321252)
(2^10,0.001662367)
(2^11,0.000831619)
(2^12,0.000415919)
};
\end{axis}
\end{tikzpicture}
\caption{Plots of the relative error in logarithmic scale.}%
\label{alphachange}%
\end{figure}

\begin{figure}[ptb]
\begin{center}
\begin{tikzpicture}[scale=0.7, transform shape]
\begin{axis}[
xlabel={$\alpha$},ylabel={$\Vert f-\tilde{f}_S\Vert_{H^{-\alpha/2}}/ \Vert f \Vert_{H^{-\alpha/2}}$},ymode=log, legend entries={Galerkin $S_\Theta^1$, Galerkin $S_\Theta^0$}, legend style={at={(1/3,1/6)},anchor=south west}]
\addplot coordinates {
(0.1,5.44E-05)
(0.2,1.29E-04)
(0.3,2.26E-04)
(0.4,3.48E-04)
(0.5,4.99E-04)
(0.6,6.80E-04)
(0.7,8.94E-04)
(0.8,1.14E-03)
(0.9,1.42E-03)
};
\end{axis}
\end{tikzpicture}
\end{center}
\caption{Plot of the relative error for fixed $N=2^{10}$.}%
\label{alphataki}%
\end{figure}
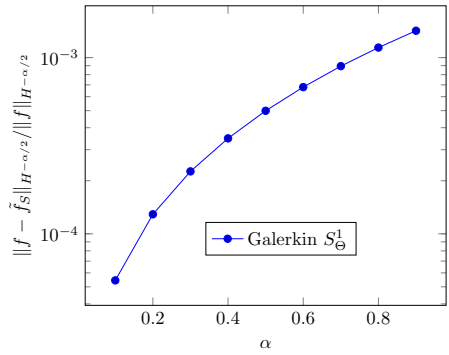

For a fixed partition size $N=2^{10}$, we display the relative error depending
on the order of singularity $\alpha$ and see the influence in Figure
\ref{alphataki}.

\section{Conclusion\label{SecConcl}}

Abel-type integral equations and its approximation have a long history and
many researchers have worked on it. This paper generalizes the idea of
Eggermont's approach in \cite{eggermont88} and develops an error analysis for
Galerkin discretizations in energy norms for generalized Abel-type kernel
functions. We propose Gauss-type quadrature methods for the approximation of
the entries of the system matrix and of the right-hand side. The local and
global error analysis allows to keep the number of quadrature nodes fairly
small. In the numerical experiments, the Galerkin method based on piecewise
polynomials of constant and first order are compared and the dependency of
error on the singularity is studied.
\bibliographystyle{acm}
\bibliography{bibli}

\end{document}